\documentclass[11pt,reqno]{amsart}
\usepackage{fullpage}
\usepackage{xcolor}
\usepackage[T1]{fontenc}
\usepackage{amsfonts}
\usepackage[utf8]{inputenc}
\usepackage{comment}
\usepackage{mparhack}
\usepackage{amsmath,amssymb,amsthm,mathrsfs,eucal}
\usepackage{mathtools}
\usepackage{booktabs}
\usepackage{graphicx,subfig}
\usepackage{wrapfig}
\usepackage[bookmarks=true,colorlinks=true]{hyperref}
\usepackage{bm}     

\usepackage{cancel}
\usepackage{soul}

\usepackage{dsfont}
\usepackage{esint}

\newtheorem{defn}{Definition}[section]
\newtheorem{thm}{Theorem}[section]
\newtheorem{prop}{Proposition}[section]
\newtheorem{lem}{Lemma}[section]
\newtheorem{cor}{Corollary}[section]
\newtheorem{ex}{Example}[section]

\DeclarePairedDelimiter{\abs}{\lvert}{\rvert}
\DeclarePairedDelimiter{\norm}{\lVert}{\rVert}

\newcommand{\mpr}{{\mathcal{P}(\R)}}

\newcommand{\mptr}{{\mathcal{P}_2(\R)}}

\newcommand{\R}{\mathbb{R}}

\newcommand{\Weps}{W_\varepsilon}
\newcommand{\Seps}{S_\varepsilon}
\newcommand{\Neps}{N_\varepsilon}

\newcommand{\sign}{\mathrm{sign}}

\newcommand{\rhoen}{\rho_\varepsilon^N}

\def\XXint#1#2#3{{\setbox0=\hbox{$#1{#2#3}{\int}$}
  \vcenter{\hbox{$#2#3$}}\kern-.5\wd0}}

\allowdisplaybreaks

\usepackage{soul}

\title{The approximation of the quadratic porous medium equation via nonlocal interacting particles subject to repulsive Morse potential}

\author{M. Di Francesco, V. Iorio, \and M. Schmidtchen}
\date{January, 2024}

\makeatletter
\@namedef{subjclassname@2020}{\textup{2020} Mathematics Subject Classification}
\makeatother

\begin{document}
\address{Marco Di Francesco - DISIM - Department of Information Engineering, Computer Science and Mathematics, University of L'Aquila, Via Vetoio 1 (Coppito)
67100 L'Aquila (AQ) - Italy}
\email{marco.difrancesco@univaq.it}

\address{Valeria Iorio - DISIM - Department of Information Engineering, Computer Science and Mathematics, University of L'Aquila, Via Vetoio 1 (Coppito)
67100 L'Aquila (AQ) - Italy}
\email{valeria.iorio1@univaq.it}

\address{Markus Schmidtchen - Institute of Scientific Computing, Technische Universitat Dresden, Zellescher Weg
25, 01217 Dresden, Germany}
\email{markus.schmidtchen@tu-dresden.de}


\maketitle

\begin{abstract}
    We propose a deterministic particle method for a one-dimensional nonlocal equation with interactions through the repulsive Morse potential. We show that the particle method converges as the number of particles goes to infinity towards weak measure solutions to the nonlocal equation. Such a results is proven under the assumption of initial data in the space of probability measures with finite second moment. In particular, our method is able to capture a measure-to-$L^\infty$ smoothing effect of the limit equation. Moreover, as the Morse potential is rescaled to approach a Dirac delta, corresponding to strongly localised repulsive interactions, the scheme becomes a particle approximation for the quadratic porous medium equation. We show that in the joint limit (localised repulsion and increasing number of particles) the reconstructed density converges to a weak solution of the porous medium equation. The strategy relies on various estimates performed at the particle level, including $L^p$ estimates and an entropy dissipation estimate, which benefit from the particular structure of our particle scheme and from the absolutely continuous reconstruction of the density from the particle locations.
\end{abstract}

\tableofcontents

\section{Introduction}

The approximation of solutions to nonlinear diffusion equations of the form 
\begin{equation}\label{eq:PME}
    \partial_t \rho = \Delta P(\rho)
\end{equation}
via a suitable combination/rearrangement of finitely many \emph{Lagrangian particles} is of interest in many contexts of applied mathematics such as population biology, micro biology, crowd management, and sampling methods. We mention here three main motivations:
\begin{itemize}
    \item [(i)] \textbf{Numerical approximation}. A successful attempt of said procedure provides a reduction of the complexity of \eqref{eq:PME} from `time dependent curves in an infinite dimensional space' to `time dependent curves in a finite dimensional space'. More precisely, one gets a \emph{particle scheme} approximation for \eqref{eq:PME} as a result.
     \item [(ii)] \textbf{Structure preserving scheme}. Partly related to the previous point, ideally the particle scheme should preserve to the greatest extent possible the main mathematical features of \eqref{eq:PME}. In the present work, we are particularly interested in the \emph{variational}, or \emph{gradient flow} structure of \eqref{eq:PME}, as well as in $L^p$ contractivity and in the $L^1$-$L^\infty$ smoothing effect.
    \item [(iii)] \textbf{Validation of continuum modelling}. In applied contexts where the unknown of \eqref{eq:PME} is the density of individuals in a population of individuals (for example in cell biology or in social sciences), a successful approximation of \eqref{eq:PME} via finitely many \emph{moving agents} validates the use of the continuum PDE \eqref{eq:PME}  in an applied framework which is intrinsically \emph{discrete}. As a result, well known qualitative results about \eqref{eq:PME} describe, with good approximation, the behavior of a large number of individuals in the discrete modelling framework.
\end{itemize}
A case which has attracted particular interest over the last decades, which is at the core of the present study, is quadratic case $P(\rho)=\rho^2/2$. Then,  \eqref{eq:PME} may be formally written as
\begin{equation}\label{eq:QPME}
    \partial_t \rho = \mathrm{div}(\rho \nabla \rho)\,.
\end{equation}
We refer at this stage to \cite{vazquez_book_1,vazquez_book_2} as comprehensive references for the mathematical theory of \eqref{eq:PME} and \eqref{eq:QPME} in particular. The first attempts at approximating diffusion equations via deterministic Lagrangian particles dates back to the pioneering paper by Russo \cite{russo} for the linear case $P(\rho)=\rho$. The nonlinear case in one space dimension is part of the results in  \cite{gosse_toscani, matthes_osberger}, see also the recent \cite{daneri_radici_runa}. In the aforementioned results, the nonlinear diffusion is represented at a discrete level via a nearest neighbor interaction, thus reproducing at a microscopic scale the `locality' in space of a diffusion mechanism.

Partly related to the Lagrangian particle approximation of \eqref{eq:QPME}, a rigorous link has been established at various stages between the quadratic porous medium equation \eqref{eq:QPME} and the \emph{nonlocal interaction equation}
\begin{equation}\label{eq:nonlocal_PME}
    \partial_t \rho = \mathrm{div}\left(\rho \nabla U_\varepsilon\ast\rho\right),
\end{equation}
where the nonlocal \emph{interaction kernel} $U_\varepsilon:\R^d\rightarrow \R$ is of the form
\[U_\varepsilon(x)=\varepsilon^{-d}U(x/\varepsilon)\,,\]
and $U$ is a fixed, nonnegative, radially decreasing kernel with unit mass, so that $U_\varepsilon$ approaches a Dirac delta and the convolution $U_\varepsilon\ast\rho$ approximates $\rho$ for small $\varepsilon$ (in a suitable distributional sense). The radially decreasing nature of $U$, often referred to as \emph{repulsiveness}, produces a diffusion-like behavior in the evolution of $\rho$ in \eqref{eq:nonlocal_PME}, with a finite speed of propagation. This suggests \eqref{eq:nonlocal_PME} to be a potentially good approximation of the nonlinear diffusion equation \eqref{eq:QPME} which shares this property.

The nonlocal approximation \eqref{eq:nonlocal_PME} has important repercussions on the particle approximation of \eqref{eq:QPME}. Indeed, an elementary calculation shows, at least when $U\in C^1(\R^d)$, that the time-depending empirical measure
\begin{equation}\label{eq:empirical_intro}
   \mu^N(t)=\frac{1}{N}\sum_{i=1}^N\delta_{x_i(t)} 
\end{equation}
solves \eqref{eq:nonlocal_PME} in the distributional sense provided the \emph{moving particles} $x_1(t),\ldots,x_N(t)$ solve the system of ODEs
\begin{equation}\label{eq:nonlocal_particles}
    \dot{x}_i(t)=-\frac{1}{N}\sum_{k=1}^N \nabla U_\varepsilon(x_i(t)-x_k(t))\,,\qquad i=1,\ldots,N\,.
\end{equation}
This very simple remark makes the nonlocal approximation of \eqref{eq:QPME} via \eqref{eq:nonlocal_PME} relevant as a particle approximation of \eqref{eq:QPME}. In contrast to the aforementioned nearest-neighbor interaction schemes, in \eqref{eq:nonlocal_particles} every particle may interact with any other particle at a distance within the support of the kernel $U_\varepsilon$. Therefore, \eqref{eq:nonlocal_particles} is often referred to as \emph{nonlocal particle approximation} of \eqref{eq:QPME}.

The approximation of \eqref{eq:QPME} via \eqref{eq:nonlocal_PME} 
as $\varepsilon\searrow 0$ and the nonlocal particle approximation of \eqref{eq:QPME} via \eqref{eq:nonlocal_particles} are therefore closely related. Two limiting scales are involved in it: the number of particles $N$, which goes to infinity, and the parameter $\varepsilon>0$, which plays the role of an \emph{interaction range} and is typically assumed to be small in order to approximate \eqref{eq:QPME}. The seminal paper \cite{oleschlaeger} by Karl Oleschlaeger first shed a light on the joint limit $N\rightarrow +\infty$ and $\varepsilon \searrow 0$ in the case in which the ODEs \eqref{eq:nonlocal_particles} feature an additional stochastic term. We deal here with the purely \emph{deterministic} problem \eqref{eq:nonlocal_particles}, which is, per se, particularly challenging when it comes to studying both the many particle limit and the small interaction range limit. A first important result on the deterministic case dates back to a paper by Lions and Mas-Gallic \cite{lions_masgallic}. More recent results frame this problem into a wider class of \emph{blob methods} for nonlinear diffusion equations, see \cite{carrillo_craig_patacchini, craig2023blob, craig2023nonlocal}. Such a problem may also be framed at a variational level, in terms of the underlying \emph{functionals}. The Wasserstein gradient flow theory \cite{AGS,JKO}  provides a natural framework in this sense. Formally, \eqref{eq:QPME} is the Wasserstein gradient flow of the energy
\[E[\rho]=\frac{1}{2}\int \rho^2 dx\]
whereas \eqref{eq:nonlocal_PME} is the gradient flow of
\[E_\varepsilon[\rho]=\frac{1}{2}\int \rho U_\varepsilon\ast \rho dx\,.\]
Some results \cite{scardia,van_meurs} use a stability approach in the framework of Gamma convergence for specific functionals arising in the modelling of dislocations, in the spirit of \cite{sandier_serfaty}. 

Closely related to the nonlocal particle approximation of \eqref{eq:QPME}, to be somehow considered as an intermediate stage between \eqref{eq:nonlocal_particles} and \eqref{eq:QPME}, is the $N\rightarrow +\infty$ limit of \eqref{eq:nonlocal_particles} for \emph{fixed} $\varepsilon>0$, towards the nonlocal PDE \eqref{eq:nonlocal_PME}. In the seminal paper \cite{dobrusin}, the case of smooth interaction potential was covered via atomisation of an absolutely continuous initial datum for \eqref{eq:nonlocal_PME} via finitely many particles, which then evolve via \eqref{eq:nonlocal_particles}. More recently, the Wasserstein gradient flow theory \cite{AGS} provided a natural framework for measure solutions to \eqref{eq:nonlocal_PME} under suitable assumption on the kernel $U_\varepsilon$. The many particle limit of \eqref{eq:nonlocal_particles} towards \eqref{eq:nonlocal_PME} may then be seen as a stability result with respect to initial data, which holds in case the interaction potential is convex up to a quadratic perturbation, see also \cite{CDFLS}. More singular kernels are considered in \cite{carrillo_choi_hauray}, with singularities no worse than Newtonian. The one-dimensional Newtonian case was considered in \cite{bonaschi_thesis,BCDP}. These results often make use of the method of characteristic, which was also used in  \cite{bertozzi_brandman,bertozzi_laurent} to prove global existence and uniqueness in a suitable $L^p$ framework in the case of blowing-up singular repulsive kernels in arbitrary dimensions. 

Given the multiple limiting parameters, let us fix some terminology at this stage. From now on we shall refer to
\begin{itemize}
    \item [(i)] the $N\rightarrow +\infty$ limit of \eqref{eq:nonlocal_particles} towards \eqref{eq:nonlocal_PME}, with fixed $\varepsilon>0$, as the \emph{nonlocal many particle} limit;
    \item [(ii)] the $\varepsilon\searrow 0$ limit of \eqref{eq:nonlocal_PME} towards \eqref{eq:QPME} as the \emph{small interaction range} limit;
    \item [(iii)] the joint $N\rightarrow +\infty$ and $\varepsilon\searrow 0$ limit of \eqref{eq:nonlocal_particles} towards \eqref{eq:QPME} as the \emph{many particle / small interaction range} limit;
\end{itemize}

A major issue in both the small interaction range limit (ii) and the join limit (iii) is to relax as much as possible the regularity assumptions on the kernel $U_\varepsilon$, see the recent \cite{burger_esposito}. The already mentioned \cite{van_meurs} uses a Gamma convergence approach to frame this problem, which turns out to be successful for some specific singular potentials too.

A particularly interesting case is when the interaction kernel $U_\varepsilon$ has a \emph{discontinuous gradient at the origin}, for example 
\begin{equation}\label{eq:morse_intro}
U_\varepsilon(x)=\varepsilon^{-1}W(x/\varepsilon)\,,\qquad W(x)=\frac{1}{2}e^{-|x|}\,.
\end{equation}
The above kernel is known as the \emph{Morse potential} and is of interest in many applications, see \cite{chuang,levine,mogilner,topaz}.

Apart from the practical implications in the applications, using the repulsive Morse potential (as well as any radially decreasing potential lacking $C^1$ regularity at the origin) somehow obstructs an empirical measure solution such as \eqref{eq:empirical_intro} to be a solution to \eqref{eq:nonlocal_PME}. Roughly said, a discontinuous gradient of $U$ prevents the \emph{particles-remain-particles} effect, and suggests a \emph{measure-to-$L^p$ smoothing effect} to take place. It is well known that the porous medium equation \eqref{eq:PME} does satisfy an $L^1-L^\infty$ smoothing effect, see e.g. \cite{veron}. The most prominent example occurs with a Dirac delta initial condition, which is known to evolve through \eqref{eq:PME} in the so called self-similar Barenblatt profile, which belongs to $L^\infty$ for all positive times.

The only chance a small interaction range limit result can catch such a smoothing effect is if the interaction kernel forces initial point particles to immediately \emph{diffuse} and regularise for positive times into an absolutely continuous density. Moreover, we stress that with such a class of singular potentials and due to the above mentioned smoothing effect, the two problems (ii) and (iii) are distinct, since the solution to \eqref{eq:nonlocal_PME} with an empirical measure initial datum is \emph{not} expected to be an empirical measure, and therefore \eqref{eq:nonlocal_PME} itself cannot be considered as a particle approximation of \eqref{eq:QPME}.

A relatively simple case in which a smoothing effect of that sort was analised in detail for a repulsive nonlocal interaction equation is contained in the paper \cite{BCDP}, which dealt with the one-dimensional case and the Newtonian kernel $U(x)=-|x|$, see also \cite{CDEFS,DiFEF} for the case of systems with two species. In that particular case, \cite{BCDP} shows that, for a given initial condition in a suitable space of probability measures which includes combinations of Dirac deltas, the solution $\rho(\cdot,t)$ becomes $L^\infty$ for positive times. Such a result is proven by exploiting a parallel between measure solutions $\rho$ to the nonlocal equation 
\[\partial_t \rho +\partial_x (\rho \partial_x |\cdot|\ast\rho)=0\]
in the $2$-Wasserstein gradient flow sense of \cite{AGS} on the one side
and entropy solutions in the Kruzkov sense to a suitable conservation law for the cumulative distribution function of $\rho$ on the other side. In particular, the measure-to-$L^\infty$ smoothing effect in the $\rho$ variable corresponds to the evolution of initial nonphysical shocks to rarefaction waves for the corresponding nonlinear conservation law. Clearly, an $\varepsilon$-scaled version of the Newtonian kernel $U(x)=-|x|$ is no use in our small interaction range limit (ii) since $U$ is not integrable at infinity in this case. Due to its discontinuous derivative at $x=0$ and to its integrability at infinity, the scaled Morse potential \eqref{eq:morse_intro} is a natural candidate to perform the small interaction range limit (ii) in such a way to capture some sort of measure-to-$L^\infty$ smoothing effect. 

\medskip
Based on the above motivations, in this paper we shall focus on  the nonlocal PDE
\begin{equation}\label{eq:main_intro}
  \partial_t \rho = (\rho \Weps'\ast \rho)_x
\end{equation}
with $\Weps$ being the scaled Morse potential
\begin{align}
  & \Weps(x)=\varepsilon^{-1}W(\varepsilon^{-1}x)\label{eq:morse_Scaled}\\
  & W(x)=\frac{1}{2}e^{-|x|}\,.\label{eq:morse}
\end{align}
Recall $W$ is nonnegative with unit mass.

Formally, as $\varepsilon\searrow 0$, $\Weps\rightarrow \delta_0$ in the weak measure sense and \eqref{eq:main_intro} behaves like the quadratic porous medium equation
\begin{equation}\label{eq:main_PME}
  \partial_t\rho = \partial_x (\rho \partial_x \rho) = \frac{1}{2}\partial_{xx}^2(\rho^2)\,.
\end{equation}
Based on the above considerations and by similarity with the approach of \cite{BCDP} and \cite{CDFLS}, a natural candidate as many-particle scheme for \eqref{eq:main_intro} is
\begin{equation}\label{eq:wrong_scheme}
    \dot{x}_i(t)=-\frac{1}{N}\sum_{k\neq i, k=1}^N \Weps'(x_i-x_k)\,.
\end{equation}
However, such a scheme has the disadvantage of invoking the gradient of $W_\varepsilon$ which does not exists at zero, a problem which would be artificially solved by eliminating the $k=i$ index in the sum in \eqref{eq:wrong_scheme}. A more reasonable approach seems to be the one recently proposed in the context of particle systems with nonlinear mobilities in \cite{radici_stra}, in which derivatives are replaced by proper difference quotients. 

Based on the above remarks, we shall work with the particle system
\begin{equation}
    \label{eq:main_particles}
    \dot{x}_i=\frac{1}{N}\sum_{k=0}^{N-1} \frac{\Weps(x_{k+1}-x_i)-\Weps(x_k-x_i)}{x_{k+1}-x_k}\,,
\end{equation}
with $i=0,\ldots,N$. The suitability of this scheme will be clear in particular in Section \ref{sec:scheme} when deriving uniform estimates. In particular, a crucial $L^p$ contractivity property will be proven for this scheme. Moreover, this scheme may be re-written in a very convenient way, see \eqref{eq:scheme_rewritten} below, which greatly simplifies the consistency of the scheme using a piecewise constant reconstruction of the density, see \eqref{eq:discrete_density} below. Indeed, we will not base our limit procedure on the empirical measure reconstruction of the particles \eqref{eq:empirical_intro}, because we need to prove $L^p$ estimates for the approximating measure.

\medskip
The main goals achieved in this paper are:
\begin{enumerate}
    \item \textbf{Many-particle limit (i) with measure initial data and smoothing effect.} We first prove that for fixed $\varepsilon>0$ the unique gradient flow solution $\rho$ to \eqref{eq:main_intro} (in the \cite{AGS} sense) is approximated for large $N$ and in a suitable weak sense by a density reconstruction of the $N$ particle system \eqref{eq:main_particles}. Such a limit if performed \emph{with general measure initial data}. The result uses crucially that the $N$-particle density reconstruction satisfies a discrete version of a measure-to-$L^\infty$ smoothing effect, which allows to detect the unique gradient flow solution in the limit. The result is stated in Theorem \ref{thm:convergence_measure} below, in which the convergence is achieved in the weak $^*$ measure topology, see also Theorem \ref{thm:convergence1} dealing with the case of $L^p$ initial data, in which the convergence holds in the weak $L^p$ sense.
    \item \textbf{Joint many particle / small interaction range limit (iii) for $L^1 \cap L^\infty$ initial data.} We then tackle the joint limit (iii) for \eqref{eq:main_particles} and prove that, under a suitable constraint on $N$ and $\varepsilon$, the solution to the quadratic porous medium equation \eqref{eq:QPME} is approximated by the density reconstruction of \eqref{eq:main_particles} in the strong $L^2$ sense for $L^1 \cap L^\infty$ initial data. The result is stated in Theorem \ref{thm:main} below.
\end{enumerate}

The second result above has also potential repercussions on the problem (ii). However, said problem may be also solved via the technique used in the recent \cite{david} for a two species system. For simplicity, we make the choice of not dealing with the problem (ii) in this paper, which is a continuum-to-continuum limit, whereas here we shall only deal with particles-to-continuum limits.

The proof of our main results requires the introduction of several methodologies. First of all, to the best of our knowledge, we do not know of any detailed analysis of the existence and uniqueness theory for \eqref{eq:main_intro} in the Wasserstein space of probability measures. While similar cases are considered in \cite{BCDP,carrillo_ferreira_precioso}, we chose to provide a detailed existence and uniqueness theory for the sake of completeness. Moreover, we need to perform an in-depth analysis of the particle scheme \eqref{eq:main_particles} when particles are \emph{initially overlapping} requiring a review of classical gradient flow theory taken from \cite{brezis}. This is needed in order to prove the many particle limit for general initial data including, possibly singular, probability measures. Several estimates are then obtained at the level of the particle scheme, which are useful in both limiting results (1) and (2) above.

The paper is structured as follows.
\begin{itemize}
    \item In Section \ref{sec:preliminaries}, we introduce tools from optimal transport, Wasserstein distances, one-dimensional continuity equations, and gradient flows. We also introduce the main concepts of solutions for the continuum PDEs \eqref{eq:main_intro} and \eqref{eq:main_PME}. Moreover, in Theorem \ref{thm:existence_and_uniqueness_GF}, we provide the main existence and uniqueness result for gradient flow solutions to \eqref{eq:main_intro} with initial data in the space of probability measures with finite second moment.
    \item In Section \ref{sec:scheme}, we introduce and solve the discrete particle scheme \eqref{eq:main_particles}, first in the case of initially non-overlapping particles, then in the case of possibly overlapping ones. The main existence and uniqueness result is contained in Theorem \ref{thm:finite_dimensional_GF}. Other important results of this section are the smoothing effect estimate in Proposition \ref{prop:discrete_smoothing}, the main uniform estimate in Lemma \ref{lem:uniform_estimates}, and the error estimate in Proposition \ref{prop:main_moment_estimate}.
    \item In Section \ref{sec:many_particle_limit}, we analyse the many particle limit for fixed $\varepsilon>0$. We define the scheme in two separate cases ($L^p$ initial data and measure initial data) and prove the two convergence results in Theorems \ref{thm:convergence1} and \ref{thm:convergence_measure}.
    \item In Section \ref{sec:joint_limit}, we analyse the joint many particles / small interaction range limit and prove our main convergence result in Theorem \ref{thm:main}.
\end{itemize}

\section{Preliminaries}\label{sec:preliminaries}

In this section we introduce the main concepts of solutions and provide the existence and uniqueness theory for \eqref{eq:main_intro}.

\subsection{Morse potential}
We recall some trivial properties of the Morse potential \eqref{eq:morse_Scaled}-\eqref{eq:morse}. $W:\R\rightarrow \R_{\geq0}$ is even. Moreover, for all $\varepsilon>0$,
\[\int_\R W(x) \, dx = \int_\R W_\varepsilon(x) \,dx = 1\,.\]
A standard computation shows that $W_\varepsilon$ is the fundamental solution of the one-dimensional elliptic operator $-\varepsilon^2 \frac{d^2}{d x^2} +\mathbb{I}$. More precisely, 
\begin{equation}\label{eq:morse_elliptic}
    -\varepsilon^2 \Weps'' + \Weps = \delta\qquad \hbox{in $\mathcal{D}'(\R)$}\,,
\end{equation}
where $\delta$ is the Dirac delta distribution centred at $0$.

\subsection{Weak solutions of the nonlocal equation}

The equation obtained in the limit $\varepsilon\searrow 0$ is the quadratic porous medium equation \eqref{eq:QPME}, which has been extensively studied in the past fifty years or so. We refer to the monographs by J. L. Vazquez \cite{vazquez_book_1,vazquez_book_2} for an extensive coverage of the mathematical theory of nonlinear diffusion equations. As far as we are concerned, we only need to define our concept of solution for \eqref{eq:QPME}, which is formulated in a weak sense as it is well known that nonlinear diffusion equations \eqref{eq:PME} with superlinear nonlinearities $P$ feature solutions with compact support and possible lack of $C^1$ regularity at the boundary of the support. For simplicity, we shall deal with the case of non-negative solutions $\rho$ which have (constant) finite mass and belong to $L^\infty([0,+\infty);\, L^1(\R)\cap L^\infty(\R))$.

\begin{defn}[Weak solutions to \eqref{eq:QPME}]\label{def:weak_PME}
    Let $\rho_0 \in L^1(\R)\cap L^\infty(\R)$ be a.e.\ non-negative. An a.e.\ non-negative measurable function $\rho \in L^\infty([0,+\infty);\, L^1(\R)\cap H^1(\R))$ is called a \emph{weak solution} to \eqref{eq:QPME} if, for all $\varphi \in C^{1}_c(\R\times [0,+\infty))$ the following identity holds
    \begin{equation}\label{eq:weak_PME}
        \int_0^{+\infty}\int_\R \left[\rho(x,t)\varphi_t(x,t) - \rho(x,t)\rho_x(x,t)\varphi_{x}(x,t)\right] dx dt + \int_\R \rho_0(x) \varphi(x,0) dx = 0
    \end{equation}
\end{defn}

We recall that the above concept of solution is well  posed. More precisely, a unique solution $\rho$ to \eqref{eq:QPME} in the sense of Definition \ref{def:weak_PME} can be found, which is globally defined on $\R\times [0,+\infty)$, see \cite{vazquez_book_1}.

As for the nonlocal interaction equation \eqref{eq:main_intro}, a concept of weak solution in the sense of measures could be easily formulated by following similar cases treated in \cite{CDFLS,bonaschi_thesis}. In these works, the discontinuity of $W'_\varepsilon$ at zero is bypassed via a simple symmetrisation technique. However, since we shall deal with solutions which are in $L^\infty(\R)$ for positive times, we can formulate our concept of weak solutions for \eqref{eq:main_intro} as in the next definition, in which the initial condition is considered in a measure sense. To this end, let $\mpr$ denote the space of probability measures on $\R$.

\begin{defn}[Weak solutions to \eqref{eq:main_intro}]\label{def:weak_nonlocal}
    Let $\mu_0 \in \mpr$. An a.e.\ non-negative measurable function $\rho\in L^\infty([0,+\infty); L^1(\R))\cap L^\infty_{\mathrm{loc}}((0,+\infty);L^\infty(\R))$ is called a \emph{weak solution} to \eqref{eq:main_intro} if for all $\varphi \in C^{1}_c(\R\times (0,+\infty))$ the following identity holds
    \begin{equation}\label{eq:weak_nonlocal}
        \int_0^{+\infty}\int_\R \left[\rho(x,t)\varphi_t(x,t) - \rho(x,t) W'_\varepsilon\ast \rho(x,t) \partial_x \varphi(x,t)\right] dx dt +\int \varphi(x,0) d\mu_0(x)= 0\,.
    \end{equation}
\end{defn}

The existence of weak solutions in the sense of Definition \ref{def:weak_nonlocal} can be easily obtained for initial data in some $L^p$ spaces with $p>1$, for example via a vanishing viscosity approximation scheme, or via mollification of the nonlocal kernel $W_\varepsilon$. As \eqref{eq:main_intro} features a smoothing effect from $\mpr$ to $L^\infty (\R)$ (as we will prove later on in the paper), existence of solutions in the sense of Definition \ref{def:weak_nonlocal} can be obtained also in the case of initial data in $\mpr$ via an approximation argument. Uniqueness is a more delicate issue. In order to obtain a concept of solution for \eqref{eq:main_intro} yielding uniqueness we need to adopt a different viewpoint which relies on the theory of \emph{gradient flows} on the $2$-Wasserstein space, which is the goal of the next sections. The bottom line of this process is that, specifically due to the $L^\infty_{\mathrm{loc}}$ regularity that we are encoding in Definition \ref{def:weak_nonlocal}, this concept of solution implies also uniqueness, see Theorem \ref{thm:existence_and_uniqueness_GF} below.

\subsection{One-dimensional Wasserstein distances}\label{subsec:one_d_wass}

We recall that $\mpr$ denotes the space of probability measures on $\R$. Let 
\[\mptr=\left\{\mu\in \mpr\,:\,\, \int_\R \abs{x}^2 d\mu(x)<+\infty\right\}\,.\]
It is well known (see e.g.\ \cite{villani}) that 
$\mptr$ is a complete metric space when equipped with the $2$-Wasserstein distance
\[
d_2(\mu_1,\mu_2)=\left[\inf_{\gamma \in \Gamma(\mu_1,\mu_2)}\iint \abs{x-y}^2 d\gamma(x,y)\right]^{1/2},
\]
where
\[\Gamma(\mu_1,\mu_2)=\left\{\gamma\in \mathcal{P}(\R\times\R)\,:\,\, \iint_{\R\times\R} \varphi(x_i)d\gamma(x_1,x_2)=\int_\R \varphi(x_i)d\mu_i(x_i)\,,\,\,\hbox{for $i=1,2$}\right\}\]
is the set of transport plans between  $\mu_1$ and $\mu_2$.

In the one-dimensional case, for $\mu\in \mptr$ we define the \emph{cumulative distribution function} $F_\mu:\R\rightarrow [0,1]$ as
\[F_\mu(x)=\mu((-\infty,x])\]
and its \emph{pseudo-inverse} (or \emph{quantile function}) $X_\mu:[0,1]\rightarrow \R$ is defined as
\[
    X_\mu(z)=\inf\left\{x\in \R\,:\,\, F_\mu(x)\geq z\right\}.
\]
For every $\mu\in \mathcal{P}_2(\R)$, the function $F_\mu:\R\rightarrow [0,1]$ is non-decreasing and right-continuous. Moreover, $X_\mu:[0,1]\rightarrow \R$ is non-decreasing as well. Therefore, $X_\mu$ has finite total variation on compact sets, which implies it has a right-continuous representative. We also recall that
\[\mu=({X_\mu}) _{\#} \mathcal{L}^1_{[0,1]}\,,\]
where $\mathcal{L}^1_{[0,1]}$ is the one-dimensional Lebesgue measure on $[0,1]$. The above is referred to $\mu$ being the \emph{push-forward} measure of $\mathcal{L}^1_{[0,1]}$ through the map $X_\mu$, which is equivalently characterised by
\begin{equation}\label{eq:change_of_variable}
    \int_\R\varphi(x)d\mu(x) = \int_0^1\varphi(X_\mu(z))dz\,,
\end{equation}
for all $\mu$-measurable functions $\varphi:\R\rightarrow\R$ such that $\frac{|\varphi(x)|}{1+|x|^2}$ is uniformly bounded, see \cite{villani,santambrogio_book}. 
One can prove (see \cite{villani}) that for all $\mu_1,\mu_2\in \mptr$ there holds
\[d_2(\mu_1,\mu_2)=\|X_{\mu_1}-X_{\mu_2}\|_{L^2([0,1])}\,.\]
By defining 
\begin{equation}\label{eq:cone}
    \mathcal{K}=\left\{X\in L^2([0,1])\, : \, \,  \hbox{$X$ is non-decreasing}\right\}\,,
\end{equation}
we have that the map
\begin{equation}\label{eq:T_map}
    \mathcal{T}: (\mptr,d_2) \ni \mu \mapsto \mathcal{T}(\mu) = X_\mu \in (\mathcal{K},\|\cdot\|_{L^2([0,1])})
\end{equation}
is an isometric bijection. Here $ (\mathcal{K},\|\cdot\|_{L^2})$ is the standard metric subspace $\mathcal{K}$ of the normed space $L^2([0,1])$. We observe that the inverse map of $\mathcal{T}$ can be expressed as follows. Let $X\in \mathcal{K}$ and let
\begin{equation}\label{eq_F_X}
   F_X(x)=\inf\left\{z\in [0,1]\,:\,\,X(z)\geq x\right\}.
\end{equation}
Then, the corresponding measure $\mathcal{T}^{-1}(X)=\mu_X \in\mptr$ is defined as $\mu_X=\partial_x F_X$ in a distributional sense.

\subsection{One-dimensional continuity equation}

Given a one-dimensional continuity equation
\begin{equation}\label{eq:continuity}
    \partial_t \rho + \partial_x (\rho v) = 0
\end{equation}
on the whole real line with a given velocity field $v=v(x,t)$, 
standard results allow to pass from \eqref{eq:continuity} to the corresponding equation for the quantile function $X(\cdot,t)=\mathcal{T}(\rho(\cdot,t)),$ which is
\begin{equation}\label{eq:continuity_quantile}
    \partial_t X(z,t)= v(X(z,t),t)\,.
\end{equation}
This is possible under pretty general assumptions on the velocity field $v$ which we do not need to investigate here. For our purpose, we assume $v$ is globally Lipschitz with respect to $x$ uniformly in $t$, which is enough to rely on classical results.

\begin{thm}\label{thm:continuity}
    Let $T\geq 0$ and let $v\in L^\infty([0,T];W^{1,\infty}(\R))$. Assume that $\rho\in L^\infty([0,T];\mptr)$ is a weak measure solution to the continuity equation \eqref{eq:continuity} on $\R\times [0,T]$ with initial condition $\overline{\rho}\in \mptr$. Then $X(\cdot,t)=\mathcal{T}(\rho(\cdot,t))$ satisfies \eqref{eq:continuity_quantile} with initial condition $\mathcal{T}(\overline{\rho})$. Vice versa, assume $X\in L^\infty([0,T]; L^2([0,1]))$ satisfies  \eqref{eq:continuity_quantile} a.e.\ on $[0,1]\times [0,T]$ with initial condition $\overline{X}\in L^2([0,1])$. Assume further that $\overline{\rho}=\mathcal{T}^{-1}(\overline{X})\in L^\infty(\R)$.  Then, $\rho(\cdot,t)=\mathcal{T}^{-1}(X(\cdot,t))$ is a weak solution to the continuity equation \eqref{eq:continuity} with initial condition $\mathcal{T}^{-1}(\overline{X})$.
\end{thm}

The proof is postponed to the Appendix (Section \ref{appendix1}) for the reader's convenience.

\subsection{Gradient flows and their Lipschitz perturbations on Hilbert spaces}\label{subsec:hilbert}

In this subsection we adapt the theory of maximal monotone operators on Hilbert spaces developed in \cite{brezis} to our purposes. Here, $H$ is just a separable Hilbert space equipped with the inner product $(\cdot,\cdot)_H$. Later on in this paper, we shall consider the case $H=L^2([0,1])$ or the case of $H$ being a finite dimensional Euclidean space.

\begin{defn}\label{def:subdiff}
Let $\mathcal{F}:H\rightarrow (-\infty,+\infty]$ be convex, bounded from below, weakly lower semi-continuous, with domain a subset $D(\mathcal{F})\subset H$. Let $X\in D(\mathcal{F})$. The \emph{subdifferential} of $\mathcal{F}$ at the point $X$ is defined as the set $\partial_X \mathcal{F}$ with
\[
Y\in \partial_X \mathcal{F}\qquad \hbox{if and only if}\qquad
\mathcal{F}(X+P)-\mathcal{F}(X)\geq (P,Y)_{H}\qquad \hbox{for all $P\in H$}\,.
\]
\end{defn}

Now let $\mathcal{F}$ be as in Definition \ref{def:subdiff} and let $\mathcal{A}:H\rightarrow H$ be a Lipschitz continuous map. For a given $\overline{X}\in \overline{D(\mathcal{F})}$ and for $t>0$ we consider the Cauchy problem
\begin{equation}\label{eq:brezis_intro}
    \begin{cases}
       - \dot{X}(t) +\mathcal{A}(X(t)) \in \partial_{X(t)}\mathcal{F},\\
       X(0)=\overline{X}.
    \end{cases}
\end{equation}

The following theorem is a consequence of
\cite[Theorem 3.17, Remark 3.14]{brezis}.

\begin{thm}\label{thm:brezis}
    Let $\overline{X}\in  \overline{D(\mathcal{F})}$. Then, there exists a unique Lipschitz continuous curve $X(\cdot):[0,+\infty)\rightarrow H$ such that \eqref{eq:brezis_intro} is satisfied. Moreover, given two solutions $X_1$ and $X_2$ to \eqref{eq:brezis_intro} with initial conditions $\overline{X}_1$ and $\overline{X}_2$ respectively, we have the stability estimate
\begin{equation}\label{eq:brezis_stability}
    \|X_1(t)-X_2(t)\|\leq e^{Lt}\|\overline{X}_1-\overline{X}_2\|,
\end{equation}
for all $t\geq 0$, with $L$ being the Lipschitz constant of $\mathcal{A}$.
\end{thm}

As a special case of the previous theorem, we consider the functional $\mathcal{F}=\mathcal{I}_{K}$ with
\[
\mathcal{I}_{K}(X)=
\begin{cases}
    0 & \hbox{if $X\in K$},\\
    +\infty & \hbox{if $X\not\in K$},
\end{cases}
\]
where $K\subset H$ is convex and closed. In this case, $D(\mathcal{F})=K$. We also observe that $\partial_X \mathcal{I}_K=\emptyset$, in the case $X\not\in K$.

\subsection{Gradient flow solution and uniqueness of weak solutions of the nonlocal equation}\label{subsec:gradient_flow_solutions}

Formally, \eqref{eq:main_intro} can be written as the $2$-Wasserstein gradient flow of the interaction energy functional
\begin{equation}\label{eq:functional}
    \mathcal{E}_\varepsilon[\mu]= \frac{1}{2}\iint_{\R\times\R}W_\varepsilon(x-y)\, d\mu(x)\, d\mu(y)\,,
\end{equation}
on the metric space $(\mptr,d_2)$, in the sense of the theory developed in \cite{AGS} and later on in \cite{CDFLS,BCDP}. In the one-dimensional case, this concept can be formulated via the isometry $\mathcal{T}$ defined in Subsection \ref{subsec:one_d_wass} as a gradient flow on the Hilbert space $L^2([0,1])$, that is, as in \eqref{eq:brezis_intro} with $\mathcal{A}\equiv 0$ and with $\mathcal{F}=E_\varepsilon:L^2([0,1])\rightarrow [0,+\infty]$ given by
\begin{equation}\label{eq:energy_E}
    E_\varepsilon[X]=\mathcal{E}_\varepsilon[\mathcal{T}^{-1}[X]]=\frac{1}{2}\iint_{[0,1]^2} \Weps(X(z)-X(\zeta))\, d\zeta \, dz\,.
\end{equation}
The equivalence between the concept of $2$-Wasserstein gradient flow in $\mptr$ and of gradient flow on $L^2([0,1])$ in the Hilbert sense via the isometry $\mu\mapsto X_\mu$ is quite easy to show for a functional of the form
\[
    \mathcal{V}[\mu]=\frac{1}{2}\iint_{\R\times\R} V(x-y) \, d\mu(x) \, d\mu(y),
\]
with $V:\R\rightarrow \R$ even, finite second moment, convex up to a quadratic perturbation, and $C^1$, see for example \cite{Bu_DiF}. We remark in particular that in this case $\partial_X \mathcal{V}$ is a singleton and
\[\partial_X \mathcal{V}(z)=\int_0^1 V'(X(z)-X(\zeta)) \,d\zeta\,,\qquad z\in [0,1]\,.\]
The cases $V(x)=\pm |x|$ were studied in \cite{BCDP}, where it was proven that $2$-Wasserstein gradient flows and $L^2$-gradient flows are equivalent concepts. Let us highlight that the repulsive potential $V(x)=-|x|$ is not convex, not even up to a quadratic perturbation. However, one can still prove the uniqueness of solutions in the $2$-Wasserstein gradient flow sense, see \cite{bonaschi_thesis}. The key property to reach that goal is that the suitable concept of convexity ensuring uniqueness relies just on the fact that $V$ is convex when restricted to $[0,+\infty)$, also using the fact that in this case solutions get smoothed instantaneously even for singular measures as initial conditions, see also \cite{carrillo_ferreira_precioso}.

Now, considering \eqref{eq:main_intro}, we observe that
\begin{equation}\label{eq:splitting}
   \Weps(x) = \Neps(x) + \Seps(x), 
\end{equation}
with
\begin{align}
  & \Neps(x)=\frac{1}{\varepsilon}N(x/\varepsilon)\,,\quad N(x)=\frac{1}{2}(1-|x|)\label{eq:splitting_N}, \\
  & \Seps(x)=\frac{1}{\varepsilon}S(x/\varepsilon)\,,\quad S(x)=\frac{1}{2}(e^{-|x|}-1+|x|)\,.\label{eq:splitting_S}
\end{align}
We compute
\begin{align*}
  & \Weps'(x)=-\frac{1}{2\varepsilon^2}\sign(x)+\frac{1}{\varepsilon^2}S'(x/\varepsilon)\,,
\end{align*}
and we also observe
\[
S''(x)=W(x)\,,\qquad \Seps''(x)=\frac{1}{\varepsilon^3}W(x/\varepsilon),
\]
which implies in particular that $\Seps$ is convex, $W^{2,\infty}$, and with $\|\Seps''\|_{L^\infty} = \|\Weps\|_{L^\infty} = \frac{1}{2\varepsilon^3}$. Hence, the existence and uniqueness of $2$-Wasserstein gradient flows for the functional $\mathcal{E}_\varepsilon$ is essentially a perturbation of the result in \cite{bonaschi_thesis}. Indeed, the singular potential $\Neps$ has been already considered in \cite{bonaschi_thesis}, whereas the perturbation potential $\Seps$ is smooth and convex. For convenience, we shall provide the detailed proof of this fact. However, we shall not formulate the $2$-Wasserstein gradient flows in the by now classical setting as in \cite{AGS}. We will instead formulate our concept of gradient flow \emph{in the pseudo-inverse variable}, similar to the technique in \cite{CDEFS} for a $2\times 2$ systems of nonlocal equations with Newtonian interactions. Let us highlight that the technique we use to prove our existence and uniqueness result is very similar to the one we will use in the discrete setting in Section \ref{sec:scheme} below. For this reason, we prefer to present the details of the existence and uniqueness proof.

\begin{defn}[Gradient flow solutions]\label{def:GF}
    Let $\overline{\mu}\in \mptr$. A curve $\mu_t:[0,+\infty)\rightarrow \mptr$ is called a \emph{gradient flow solution} of \eqref{eq:main_intro} with initial condition $\overline{\mu}$ if
    \[\lim_{t\searrow 0}\mu_t=\overline{\mu},\]
    in the $d_2$ sense and, given $\mathcal{T}(\mu_t):[0,+\infty)\rightarrow L^2([0,1])$ the pseudo inverse obtained from $\mu_t$, the curve $X(\cdot,t)=\mathcal{T}(\mu_t)$ solves
    \begin{equation}\label{eq:GF_hilbert}
        -\dot{X}(\cdot,t) \in \partial_{X(\cdot,t)} \left(E_\varepsilon+\mathcal{I}_\mathcal{K}\right)\,,\qquad \hbox{for all $t>0$},
    \end{equation}
    where $E_\varepsilon$ is the energy defined in \eqref{eq:energy_E} and $\mathcal{I}_\mathcal{K}$ is the indicator function of the closed convex cone $\mathcal{K}$ \eqref{eq:cone}.
\end{defn}
For future reference, we use the notation $X(\cdot,t)$ to denote the $L^2([0,1])$ evaluation of the curve $\mathcal{T}(\mu_t)$ at time $t$, whereas we denote by $X(z,t)\in \R$ the value attained by the unique right-continuous representative of $X(\cdot,t)$ at some given $z\in [0,1]$.

We now state the existence and uniqueness of gradient flow solutions to \eqref{eq:main_intro} in the next theorem, which also collects further important properties, the most notable one being that weak solution in the sense of Definition \ref{def:weak_nonlocal} are \emph{unique}.

\begin{thm}\label{thm:existence_and_uniqueness_GF}
Let $\overline{\mu}\in \mptr$. Then there exists a unique gradient flow solution $\mu_t$ to \eqref{eq:main_intro} with initial condition $\overline{\mu}$ in the sense defined in Definition \ref{def:GF}. Moreover, the following properties hold:
\begin{itemize}
    \item [(i)] Given $X(\cdot,t)=\mathcal{T}(\mu_t)$, $X(\cdot,t)$ is strictly increasing on $[0,1]$ for positive times and satisfies
\begin{equation}\label{eq:GF_Wprime}
    \dot{X}(z,t)=-\int_0^1\Weps'(X(z,t)-X(\zeta,t))d\zeta, \qquad \hbox{for a.e.\ $(z,t)\in[0,1]\times (0,+\infty)$}.
\end{equation}
\item [(ii)] For a.e. $t>0$, $\mu_t$ is absolutely continuous with respect to the one-dimensional Lebesgue measure and belongs to $L^\infty(\R)$. Moreover, the density $\rho(x,t)$ of $\mu_t$ is a weak solution to \eqref{eq:main_intro} with initial condition $\overline{\mu}$ in the sense of Definition \ref{def:weak_nonlocal}.
\item [(iii)] Given a curve $\mu_t:[0,+\infty)\rightarrow\mptr$ with $\mu_0=\overline{\mu}$ and such that $X(\cdot,t)=\mathcal{T}(\mu_t)$ is strictly increasing on $[0,1]$ for positive times and satisfies \eqref{eq:GF_Wprime}, then $\mu_t$ is the unique gradient flow solution to \eqref{eq:main_intro} with initial datum $\overline{\mu}$.
\item [(iv)] Given two gradient flow solutions $\mu_t,\nu_t$ to \eqref{eq:main_intro} with initial conditions $\overline{\mu}$ and $\overline{\nu}$ respectively, the following stability estimate holds:
\begin{equation}\label{eq:stability_main}
    d_2(\mu_t,\nu_t)\leq d_2(\overline{\mu},\overline{\nu})\,,\qquad t\geq 0\,.
\end{equation}
\item [(v)] There exists a unique weak solution to \eqref{eq:main_intro} in the sense of Definition \ref{def:weak_nonlocal} for a given initial condition $\overline{\mu}\in \mptr$.
\end{itemize}
\end{thm}

We prove Theorem \ref{thm:existence_and_uniqueness_GF} in several steps. In what follows, $\overline{X}=\mathcal{T}(\overline{\mu})$. Consistently with the decomposition \eqref{eq:splitting}-\eqref{eq:splitting_N}-\eqref{eq:splitting_S}, we set
\[\mathcal{N}_\varepsilon[X]=\frac{1}{2}\int_0^1\int_0^1 \Neps(X(z)-X(\zeta))d\zeta dz\,.\]
Since the potential $N_\varepsilon$ is not convex, the convexity of the functional $\mathcal{N}_\varepsilon$ is not immediate. However, 
we have the following lemma, the proof of which is contained in \cite{natile} and \cite{BCDP}. We provide the statement here and postpone the proof to the Appendix (Section \ref{appendix1}) for the reader's convenience.

\begin{lem}\label{lem:N_and_R_coincide}
    Let
    \[\mathcal{R}_\varepsilon[X]= \frac{1}{4\varepsilon}\left[1-\frac{2}{\varepsilon}\int_0^1(2z-1)X(z) dz\right]\,.\]
    Then, $\mathcal{N}_\varepsilon$ and $\mathcal{R}_\varepsilon$ coincide on the cone $\mathcal{K}$ defined in \eqref{eq:cone}. Consequently,
    \[\mathcal{N}_\varepsilon + \mathcal{I}_\mathcal{K} =\mathcal{R}_\varepsilon + \mathcal{I}_\mathcal{K}\,.\]
\end{lem}

As a consequence of Lemma \ref{lem:N_and_R_coincide}, we can replace \eqref{eq:GF_hilbert} by the differential inclusion
\begin{equation}\label{eq:GF_reformulated}
  -\dot{X}(\cdot,t)\in \partial_{X(\cdot,t)}\left(\mathcal{S}_\varepsilon+\mathcal{R}_\varepsilon+\mathcal{I}_{\mathcal{K}}\right)\,, 
\end{equation}
where
\[\mathcal{S}_\varepsilon[X]=\frac{1}{2}\int_0^1 \int_0^1 \Seps(X(z)-X(\zeta))d\zeta dz\,.\]
We observe that 
$\partial_{X}\mathcal{R}_\varepsilon$ is a singleton for all $X\in L^2([0,1])$ and 
\[
    \partial_{X}\mathcal{R}_\varepsilon(z)=-\frac{2z-1}{2\varepsilon^2}, \qquad \hbox{for all $X\in L^2([0,1])$}.
\]
Moreover, as a linear functional, $\mathcal{R}_\varepsilon$ is trivially convex. Since $\Seps$ is continuously differentiable, its subdifferential is a singleton given by
\[
    \partial_X\mathcal{S}_\varepsilon(z)=\int_0^1 \Seps'(X(z)-X(\zeta))d\zeta\,,\qquad z\in [0,1].
\]
We also observe that $\mathcal{S}_\varepsilon$ is convex, as $S_\varepsilon''\geq 0$. Hence, \eqref{eq:GF_reformulated} can be re-written, for $t>0$, as
\begin{equation}\label{eq:auxiliary_GF}
    \begin{cases}
        -\displaystyle{\dot{X}(\cdot,t)+ \frac{2z-1}{2\varepsilon^2} -\int_0^1 \Seps'(X(z,t)-X(\zeta,t))d\zeta} \in \partial_{X(t)}\mathcal{I}_{\mathcal{K}},  \\
        X(\cdot,0)=\overline{X}\,.
    \end{cases}
\end{equation}
The first identity in \eqref{eq:auxiliary_GF} is set on $L^2([0,1])$, therefore the evaluation on $z$ in the integral term is a slight abuse of notation. 

Since all three functionals above are convex, Theorem \ref{thm:brezis} implies the existence and uniqueness of a solution to \eqref{eq:auxiliary_GF}. 
Moreover, since $\partial_{X(\cdot,t)}\mathcal{I}_{\mathcal{K}}\neq \emptyset$ for all $t>0$, we have that $X(\cdot,t)\in \mathcal{K}$ for all $t\geq 0$.

We now prove the following Lemma.

\begin{lem}\label{lem:lemma_pseudo_increasing}
    The solution $X(\cdot,t)$ to \eqref{eq:auxiliary_GF} is strictly increasing with respect to $z\in [0,1]$ for all $t>0$, and the following estimate holds
    \begin{equation}\label{eq:slope_pseudo}
        \frac{X(z_2,t)-X(z_1,t)}{z_2-z_1}\geq \frac{1}{\varepsilon}\left(1-e^{-\frac{t}{2\varepsilon^3}}\right)\,,\qquad 0\leq z_1<z_2\leq 1\,.
    \end{equation}
\end{lem}

\begin{proof}
    Let $z_1,z_2\in [0,1]$ with $z_1<z_2$. Since $X(t)\in \mathcal{K}$ for all $t\geq 0$, the function $X(\cdot,t)$ is non-decreasing for all times. Hence,
\begin{align*}
    & \frac{d}{dt}\left(X(z_2,t)-X(z_1,t)\right)\\
    & \ = \frac{z_2-z_1}{\varepsilon^2}-\int_0^1 \left(\Seps'(X(z_2,t)-X(\zeta,t))-\Seps'(X(z_1,t)-X(\zeta,t))\right) d\zeta\\
    & \ \geq \frac{z_2-z_1}{\varepsilon^2}-\frac{1}{2\varepsilon^3}\left(X(z_2,t)-X(z_1,t)\right),
\end{align*}
and Gr\"onwall's inequality implies
\begin{align*}
    & X(z_2,t)-X(z_1,t)\geq e^{-\frac{t}{2\varepsilon^3}}(X(z_2,0)-X(z_1,0))+\frac{z_2-z_1}{\varepsilon^2}\left(1-e^{-\frac{t}{2\varepsilon^3}}\right)\\
    & \geq \frac{z_2-z_1}{\varepsilon^2}\left(1-e^{-\frac{t}{2\varepsilon^3}}\right)>0
\end{align*}
if $t>0$, that proves the result.
\end{proof}

As a consequence of Lemma \ref{lem:lemma_pseudo_increasing}, we have the following

\begin{lem}\label{lem:singular_linear}
Let $X(\cdot,t)$ be the unique solution to the auxiliary Cauchy problem \eqref{eq:auxiliary_GF}. Then, 
    for all $t>0$, the function
    \[ [0,1]\ni z\mapsto \partial_{X(\cdot,t)}\mathcal{R}_\varepsilon(z) = -\frac{2z-1}{2\varepsilon^2}\]
    satisfies
    \[\partial_{X(\cdot,t)}\mathcal{R}_\varepsilon(z) =
-\frac{1}{2\varepsilon^2}\int_0^1 \sign(X(z,t)-X(\zeta,t))d\zeta\,.
    \]
\end{lem}

\begin{proof}
   Lemma \ref{lem:lemma_pseudo_increasing} implies for $t>0$
\begin{align*}
    & \int_0^1 \sign(X(z,t)-X(\zeta,t))d\zeta = \int_0^z d\zeta -\int_z^1 d\zeta = 2z-1\,.
\end{align*}
Therefore,
\[\partial_{X(\cdot,t)}\mathcal{R}_\varepsilon(z) = -\frac{2z-1}{2\varepsilon^2} = -\frac{1}{2\varepsilon^2}\int_0^1 \sign(X(z,t)-X(\zeta,t))d\zeta\]
    which concludes the proof.
\end{proof}

\begin{proof}[Proof of Theorem \ref{thm:existence_and_uniqueness_GF}]\mbox{}\\
\noindent
\underline{Ad existence of unique gradient flow:} 
    The differential inclusion \eqref{eq:GF_hilbert} is equivalent to \eqref{eq:GF_reformulated}. Since $\mathcal{R}_\varepsilon$ is linear and $\mathcal{S}_\varepsilon$ and $\mathcal{I}_\mathcal{K}$ are convex, Theorem \ref{thm:brezis} implies the existence and uniqueness of a solution $X(\cdot,t)$ to \eqref{eq:auxiliary_GF}, and hence to \eqref{eq:GF_hilbert}.
    This proves the opening statement of the theorem. 

\noindent
\underline{Ad (i):}
    Lemma \ref{lem:singular_linear} shows the identity \eqref{eq:GF_Wprime} is satisfied for positive times, which proves (i). 

\noindent
\underline{Ad (ii):}
Moreover, Lemma \ref{lem:lemma_pseudo_increasing} shows $X(\cdot,t)$ is strictly increasing for positive times. Hence, given $\mu_t=\mathcal{T}^{-1}(X(\cdot,t))$, estimate \eqref{eq:slope_pseudo} implies $F_{\mu_t}$ is globally Lipschitz on $\R$ for each positive time $t$ and $\rho(\cdot,t)=\partial_x F_{\mu_t}\in L^\infty(\R)$. In order to prove that $\rho$ is a weak solution in the sense of Definition \ref{def:weak_nonlocal}, we start from \eqref{eq:GF_Wprime}, which can be written as
    \[\partial_t X(z,t)=v(X(z,t),t),\]
    with
    \[v(x,t)=-\int_0^1 W'_\varepsilon(x-X(\zeta,t))d\zeta.
    \]
    Hence, since $\rho(\cdot,t_0)\in L^\infty(\R)$ for an arbitrary $t_0>0$, Theorem \ref{thm:continuity} implies that $\rho(\cdot,t)$ satisfies the definition of weak solution in Definition \ref{def:weak_nonlocal} on $\R\times [t_0,+\infty)$. Indeed, due to \eqref{eq:change_of_variable} $v$ satisfies
\begin{equation}\label{eq:v_continuity}
    v(x,t)=-\int_\R W'_\varepsilon(x-y)\rho(y,t) dy\,.
\end{equation}
    Since $\rho(\cdot,t)\in L^\infty(\R)$ for $t\in [t_0,+\infty)$, we can differentiate with respect to $x$ and get
    \[v_x(\cdot,t)=-W''_\varepsilon\ast \rho(\cdot,t)\,.\]
    Hence, $v$ is uniformly bounded and with uniformly bounded Lipschitz norm on bounded intervals $(t_1,t_2)$.
 This proves (ii). 
 
\noindent
\underline{Ad (iii):}
Since \eqref{eq:GF_Wprime} is satisfied and $X(\cdot,t)$ is strictly increasing for $t>0$, Lemma \ref{lem:singular_linear} and the identity $\mathcal{R}_\varepsilon[X(\cdot,t)]=\mathcal{N}_\varepsilon[X(\cdot,t)]$ imply $X(\cdot,t)$ solves \eqref{eq:GF_hilbert} for $t>0$ and satisfies $X(\cdot,0)=\overline{X}$. This means $\mu_t$ is a gradient flow solution, proving (iii). 
 
\noindent
\underline{Ad (iv):} Statement (iv) is just a consequence of Theorem \ref{thm:brezis}. 

\noindent 
\underline{Ad (v):} To prove (v), assuming $\rho$ is as in Definition \ref{def:weak_nonlocal} (in particular with the required $L^\infty$ regularity in space), we observe that Theorem \ref{thm:continuity} implies that $X(\cdot,t)=\mathcal{T}(\rho(\cdot,t))$ satisfies \eqref{eq:continuity_quantile} with $v$ as in \eqref{eq:v_continuity} a. e. on $t\in [\delta,+\infty)$ for an arbitrary $\delta>0$, and hence a. e. on $t\in (0,+\infty)$. By choosing a test function $\varphi_\delta(x,t)=\xi(x)\psi_\delta(t)$ with $\xi\in C^1_c(\R)$ and with $\psi_\delta\in C^1_c([0,+\infty))$ with $\psi_\delta \geq 0$, $\psi(0)=1$, $\psi(t)=0$ for $t\in (\delta,+\infty)$, $\psi_\delta$ non increasing on $[0,\delta]$ and $\psi'_\delta\rightarrow -\delta_0$ in the sense of distributions, we substitute $\varphi_\delta$ in \eqref{eq:weak_nonlocal} and obtain
 \begin{align*}
     &  \int_0^{+\infty}\psi_\delta'(t)\left(\int_\R \rho(x,t)\xi(x) dx\right) dt+ \int_0^{\delta} \psi(t)\left(\int_R\rho(x,t) W'_\varepsilon\ast \rho(x,t) \xi'(x)dx\right) dt +\int\xi(x)d\mu_0(x)= 0.
 \end{align*}
We note that $\Weps'\in L^\infty(\R)$ implies $\Weps'\ast\rho\in L^\infty( \R\times [0,\delta])$, and by letting $\delta \to 0$, we obtain
\begin{align*}
     & \lim_{t\searrow 0}\int_\R \rho(x,t)\xi(x) dx = \int \xi(x) d\mu_0(x)\,,
\end{align*}
which proves that $\rho$  is a unique gradient flow solution.
\end{proof}

\section{The discrete particle model: well-posedness and estimates}\label{sec:scheme}

In this section, we introduce the particle model \eqref{eq:main_intro} with a fixed number of particles in detail. We start with the case in which particles are initially strictly ordered, and then shift our focus on the case of initially (possibly) overlapping particles.

\subsection{The particle scheme for strictly ordered particles}\label{subsec:scheme_strictly_ordered}

For $x\in \R^{N+1}$ we use the notation
\[x=(x_0,\ldots,x_N),\]
where $x_i$, $i=0,\ldots,N$, denote the positions of the particles.
For a fixed $N\in \mathbb{N}$ we define the open set of non-overlapping positions
\[O_N=\left\{x\in \R^{N+1}\,:\,\, x_i<x_{i+1}\,\,\hbox{for all $i=0,\ldots,N-1$}\right\}\,.\]
We consider an initial vector $\overline{x}\in O_N$. We then solve the dynamical system
\begin{equation}
    \label{eq:scheme}
    \dot{x}_i=\frac{1}{N}\sum_{k=0}^{N-1} \frac{\Weps(x_{k+1}-x_i)-\Weps(x_k-x_i)}{x_{k+1}-x_k}\,,
\end{equation}
with initial conditions 
\begin{equation}\label{eq:scheme_initial}
x_i(0)=\overline{x}_i\,,\qquad i=0,\ldots,N\,.
\end{equation}
Since $\Weps$ is a Lipschitz function and the differences $x_{k+1}-x_k$ are initially positive, the Cauchy problem \eqref{eq:scheme}-\eqref{eq:scheme_initial} admits a unique solution $x(t)=(x_0(t),\ldots,x_N(t))\in O_N$ with $t\in [0,t^*)$ and $t^*$ a sufficiently small (maximal) local existence time, by the classical Cauchy-Lipschitz theorem. In order to prove that $t^*$ can be extended to $+\infty$ we show that particles never collide for finite times. More precisely, we prove that on arbitrary finite time intervals $[0,T]$, the distance between any two consecutive particles remainy bounded away from zero uniformly with respect to $t\in [0,T].$ Then, a classical continuation principle implies that the solution to \eqref{eq:scheme}-\eqref{eq:scheme_initial} stays in $O_N$ for all times.

To obtain said non-collision property, we introduce the differences
\begin{equation}\label{eq:differences}
    d_i=x_{i+1}-x_i\,,\qquad i=0,\ldots,N-1
\end{equation}
and the discrete densities
\begin{equation}\label{eq:Ri}
   R_i=\frac{1}{N d_i} = \frac{1}{N(x_{i+1}-x_i)}\,,\qquad i=0,\ldots,N-1\,. 
\end{equation}
We then define the piecewise constant density
\begin{equation}\label{eq:discrete_density}
    \rho^N(x,t)=\sum_{k=0}^{N-1}R_k(t)\mathbf{1}_{[x_k(t),x_{k+1}(t))}(x),
\end{equation}
and observe that 
\begin{align*}
    \dot{x}_i 
    &=\frac{1}{N}\sum_{k=0}^{N-1} \frac{\Weps(x_{k+1}-x_i)-\Weps(x_k-x_i)}{x_{k+1}-x_k}\\
    &= \sum_{k=0}^{N-1} \int_{x_k}^{x_{k+1}} \frac1N \frac{1}{d_k} \Weps'(y - x_i) dy\\
    & = \sum_{k=0}^{N-1} \int_{x_k}^{x_{k+1}} \Weps'(y-x_i) \rho^N(y,t) dy\\
    & \ = \int_{\R}\Weps'(y-x_i)\rho^N(y,t) dy\,,
\end{align*}
which implies that the ODE system in \eqref{eq:scheme} can actually be rewritten in the very natural form
\begin{equation}\label{eq:scheme_rewritten}
\dot{x}_i(t)= -\Weps'\ast \rho^N  (x_i(t))\,.
\end{equation}
For future use, we provide here the following identity for the evolution of $d_i$
\begin{equation}\label{eq:d_i_dot}
    \dot{d}_i(t)=-\Weps'\ast \rho^N  (x_{i+1}(t),t)+\Weps'\ast \rho^N  (x_i(t),t) = -\int_{x_i(t)}^{x_{i+1}(t)}\Weps''\ast \rho^N(y,t) dy,
\end{equation}
and $R_i$
\begin{equation}\label{eq:R_i_dot}
    \dot{R}_i(t)=-\frac{\dot{d}_i(t)}{N d_i(t)^2} = R_i(t)\fint_{x_i(t)}^{x_{i+1}(t)}\Weps''\ast \rho^N(y,t) dy\,.
\end{equation}
In the following lemma we provide a basic estimate ensuring that particles do not collide at finite times, reminiscent of Lemma \ref{lem:lemma_pseudo_increasing}.

\begin{lem}[Particles do not collide]\label{lem:particles_do_not_collide}
Let $\bar{x}\in O_N$. Then, for all $t>0$ we have the estimate
\begin{equation}\label{eq:estimate_d_i}
    d_i(t)\geq d_i (0) e^{-\frac{t}{2\varepsilon^3}} + \frac{2\varepsilon}{N}  \bigg[ 1-e^{-\frac{t}{2\varepsilon^3}} \bigg].
\end{equation}
In particular, the unique local-in-time solution $x(t)$ to \eqref{eq:scheme}-\eqref{eq:scheme_initial} belongs to $O_N$ for all $t>0$ and consequently it exists for all times.
\end{lem}

\begin{proof}
Using \eqref{eq:d_i_dot}, we get
\begin{align*}
\dot{d}_i =& -\int_{x_i}^{x_{i+1}} \Weps'' \ast \rho^N (y)\,dy  \\
=&\frac{1}{\varepsilon^2} \int_{x_i}^{x_{i+1}} (\delta_0 -\Weps )\ast \rho^N (y)\,dy \\
=& \frac{1}{\varepsilon^2} \int_{x_i}^{x_{i+1}} \rho^N (y)\,dy + \frac{1}{\varepsilon^2} \int_{x_i}^{x_{i+1}} \Weps \ast \rho^N (y)\,dy \\
=& \frac{1}{\varepsilon^2} \bigg[ \frac{1}{N} - \int_{x_i}^{x_{i+1}} \Weps \ast \rho^N (y)\,dy \bigg].
\end{align*}
By Young's inequality, we have
\[ 
    \norm{\Weps \ast \rho^N}_{L^\infty} \leq \norm{\Weps}_{L^\infty} \norm{\rho ^N}_{L^1} \leq \frac{1}{2\varepsilon},
\]
and then
\begin{equation}\label{eq:di_dot}
   \dot{d}_i \geq \frac{1}{\varepsilon^2} \bigg[ \frac{1}{N} - \frac{1}{2\varepsilon} \bigg] = \frac{1}{N\varepsilon^2} - \frac{d_i}{2 \varepsilon^3}. 
\end{equation}
Integrating in time, we get
\begin{align*}
    d_i (t) \geq \ & e^{-\frac{t}{2\varepsilon^3}} \bigg[ d_i (0) + \int_0^t \frac{1}{N\varepsilon^2} e^{\frac{s}{2\varepsilon^3}}\,ds \bigg] \\
    = \ & d_i (0) e^{-\frac{t}{2\varepsilon^3}} + \frac{2\varepsilon}{N}  \bigg[ 1-e^{-\frac{t}{2\varepsilon^3}} \bigg].
\end{align*}
\end{proof}

Lemma \ref{lem:particles_do_not_collide} implies the global-in-time well-posedness of the particle scheme \eqref{eq:scheme} with initially non-overlapping positions \eqref{eq:scheme_initial}. Indeed, \eqref{eq:scheme_rewritten} shows that $\dot{x}_i(t)$ is uniformly bounded on given bounded time intervals for fixed $\varepsilon$, which, upon integrating, implies a global uniform bound for all particles $x_i(t)$ themselves. We now use estimate \eqref{eq:estimate_d_i} to provide $L^p$ estimates on $\rho^N$ which are uniform in time. Such estimates automatically translates into a uniform-in-$N$ estimates in case the initial $L^p$ norm of $\rho^N$ is uniformly bounded in $N$.

\begin{lem}[$L^p$ estimates for $\rho^N$]\label{lem:estimates_1} Let $p \in [1, \infty]$. Then, $\norm{\rho^N(t)}_{L^p}$ is bounded uniformly in time. More precisely
\begin{equation}\label{eq:Lp_estimate_1}
    \norm{\rho^N (t)}_{L^p} \leq \max \bigg\{ \norm{\rho^N (0)} _{L^p}, \frac{1}{(2\varepsilon)^{1- \frac{1}{p}}} \bigg\}\,,
\end{equation}
with the convention $1/\infty = 0$.
\end{lem}
\begin{proof}
From Lemma \ref{lem:particles_do_not_collide}, we know that
\begin{align*}
R_i (t) =  \frac{1}{N d_i(t)} \leq & \frac{1}{N d_i(0) e^{- \frac{t}{2\varepsilon^2}}+2\varepsilon ( 1- e^{-\frac{t}{2\varepsilon^3}})} \\
= & \frac{R_i (0)}{e^{-\frac{t}{2\varepsilon^2}} + 2\varepsilon R_i (0) (1-e^{-\frac{t}{2\varepsilon^2}})}.
\end{align*}
Using the fact that the denominator is a convex combination of $1$ and $2\varepsilon R_i(0)$, we get
\begin{align*}
R_i (t) \leq \ & R_i(0) \frac{1}{\min\{ 1, 2\varepsilon R_i (0) \}} \\
\leq \ & R_i(0) \max \bigg\{ 1, \frac{1}{2\varepsilon R_i (0)} \bigg\} \\
\leq \ & \max \bigg\{ R_i(0), \frac{1}{2\varepsilon} \bigg\}.
\end{align*}
This implies for $t>0$
\[ \norm{\rho^N (t)}_{L^\infty} = \max_i R_i(t) \leq \max \bigg\{ \max_i R_i (0), \frac{1}{2\varepsilon} \bigg\}, \]
that is
\[ \norm{\rho ^N (t)}_{L^\infty} \leq \max \bigg\{ \norm{ \rho^N(0)}_{L^\infty}, \frac{1}{2\varepsilon} \bigg\}. \]
Now let $p\geq 1$.
\begin{align*}
    \norm{\rho^N}_{L^p}^p & = \int_\R \abs{\rho^N (x)}^p\,dx = \sum_{i=0}^{N-1} \int_{x_i}^{x_{i+1}} R_i (t)^p\,dx \\
    & = \sum_{i=0}^{N-1} R_i (t)^p d_i = \sum_{i=0}^{N-1} \frac{1}{N} R_i (t)^{p-1}.
\end{align*}
From the previous computation, we obtain
\begin{align*}
    \norm{\rho^N}_{L^p}^p & \leq \sum_{i=0}^{N-1} \max \bigg\{ \frac{1}{N} R_i(0)^{p-1}, \frac{1}{(2\varepsilon)^{p-1}} \frac{1}{N} \bigg\} \\
    & \leq \max \bigg\{ \sum_{i=0}^{N-1} R_i(0)^{p-1} R_i(0) d_i(0), \sum_{i=0}^{N-1} \frac{1}{(2\varepsilon)^{p-1}} \frac{1}{N} \bigg\} \\
    & = \max \bigg\{ \norm{\rho^N (0)}_{L^p}^p, \frac{1}{(2\varepsilon)^{p-1}} \bigg\}.
\end{align*}
Hence, \eqref{eq:Lp_estimate_1} is proven.
\end{proof}

\subsection{Particle scheme for initially overlapping particles}\label{subsec:scheme_overlapping}

We now extend the well-posedness theory for \eqref{eq:scheme} performed in Subsection \ref{subsec:scheme_strictly_ordered} to the case in which particles are possibly overlapping at time $t=0$. Recalling the notation for $O_N\subset \R^{N+1}$ in Subsection \ref{subsec:scheme_strictly_ordered}, we set
\[
    K_N:=\overline{O_N}=\left\{x\in \R^{N+1}\,:\,\, x_i\leq x_{i+1}\,\,\hbox{for all $i=0,\ldots,N-1$}\right\}\,.
\]
The set $K_N$ is closed and convex in $\R^{N+1}$. Here $\R^{N+1}$ is equipped with the following inner product structure: given $x,y\in \R^{N+1}$, $x=(x_0,\ldots,x_N)$, $y=(y_0,\ldots,y_N)$, we set
\[
    \langle x,y\rangle_N := \frac{1}{N}\sum_{i=0}^{N} x_i y_i\,.
\]
As customary, we then set
\begin{equation}
    \label{eq:scaled-norm}
    \|x\|_N:=\sqrt{\langle x,x\rangle_N}\,.
\end{equation}

The main difficulty in dealing with initially overlapping particles resides in the fact that one might either ignore self-interactions (that is, interaction with particles located in the same position) or exploit the repulsive feature of the interaction potential and \emph{force} particles to detach. In the next example we illustrate our argument in favour of the latter option.

\begin{ex}\label{ex:two_particles}
\emph{
We motivate our strategy with the following example. Consider $N=1$ and let $\overline{x}_0=\overline{x}_1=0$. The ODE particle system in this case reads
\begin{align*}
    & \dot{x}_1(t)=\frac{-\Weps(x_1(t)-x_0(t))+\Weps(0)}{x_1(t)-x_0(t)} = -\dot{x}_0(t)\,,
\end{align*}
provided $x_0(t)<x_1(t)$. Clearly, one could also consider the case in which the two particles maintain the same initial positions, but this requires dealing with difference quotients of $\Weps$ at zero, which are in principle not defined in the limit. A reasonable way to bypass the lack of continuity at zero of $\Weps'$ is to remove all interactions with particles located at the same position as $x_i$ in the ODE \eqref{eq:main_particles}. In this particular example, this would imply
\[\dot{x}_0(t)=\dot{x}_1(t)=0\,.\]
Hence, we potentially have two distinct, meaningful solutions with the same initial condition. To single out the physically relevant one, we observe that the first solution is \emph{stable} with respect to any perturbation that preserves the mid-point between $x_0$ and $x_1$, whereas the second solution is not. Indeed, by slightly detaching the two particles initially, they will move apart as in the first solution case, which makes the second solution unstable. 
}
\end{ex}

\noindent\textbf{Relation to gradient flows with Newtonian interactions.} Based on the above example, we aim at formulating a concept of solution which allows to single out a unique solution and according to which $k$ particles initially overlapping split instantaneously into $k$ particles with distinct positions. Such a concept of solution relates to \emph{gradient flows} on the finite dimensional space $\R^{N}$. It is well known that if $F:\R\rightarrow \R$ is a convex, smooth functional growing less than quadratically at infinity, then the gradient flow of the energy
\[
    \mathcal{F}[\mu]=\frac{1}{2N^2}\sum_{i=1}^{N}\sum_{j=1}^{N}F(x_i-x_j),
\]
on $\R^{N}$ equipped with the scaled Euclidean norm $\|\cdot\|_N$ is formally given by the particle system
\begin{equation}\label{eq:scheme_GF_F}
   \dot{x}_i(t)=-\frac{1}{N}\sum_{k=1}^N F'(x_i-k_k)\,. 
\end{equation}
The choice $F(x)=\frac{1}{2}(1-|x|)$ is part of the theory treated in detail in \cite{BCDP}, where the main issue is the non-convexity of $F$, not even up to a quadratic perturbation, because of the repulsive Lipschitz singularity at $x=0$. However, \cite{BCDP} showed that if one prescribes \emph{a priori} that particles $x_1<\ldots<x_N$ are not overlapping for positive times (which is reasonable when the energy is repulsive) and are ordered via their label $i=0,\ldots,N$, one can rephrase the particle system corresponding to the potential $F(x)=-\frac{1}{2}(1-|x|)$ as
\begin{align*}
    & \dot{x}_i=\frac{1}{2N}\sum_{k=0\,,\,k\neq i}^{N-1}\sign(x_i-x_k) = \frac{1}{2N}\sum_{k=0\,,\,k\neq i}^{N-1} \sign (i-k)\\
    & \ = \frac{1}{2N}\sum_{k=0}^{i-1} 1 -\frac{1}{2N}\sum_{k=i+1}^{N-1} 1 = \frac{2i-N+1}{2N}.
\end{align*}

\noindent\textbf{Morse potential setting.} Let us now get back to the case of our interaction potential $W_\varepsilon$. The expressions in \eqref{eq:splitting}-\eqref{eq:splitting_N}-\eqref{eq:splitting_S} suggest that, similar to Subsection \ref{subsec:gradient_flow_solutions}, we can treat $W_\varepsilon$ as a smooth and convex perturbation of $N_\varepsilon$, the latter being exactly the same as the case treated in \cite{BCDP}, discussed above. However, our particle scheme \eqref{eq:scheme} has a slightly different structure than \eqref{eq:scheme_GF_F}, in that the latter uses the derivative of interaction potential $W_\varepsilon$, whereas we use difference quotients of $W_\varepsilon$ in \eqref{eq:scheme}. Nevertheless, we prove that a theory consistent with Example \ref{ex:two_particles} can also be developed for \eqref{eq:scheme} via the splitting \eqref{eq:splitting}. Here, the term $S_\varepsilon$ is treated as a smooth perturbation of $N_\varepsilon$ instead of as an additional gradient flow term. Based on \eqref{eq:splitting}-\eqref{eq:splitting_N}-\eqref{eq:splitting_S}, we recall
$S''(x)=W(x)$, 
which implies in particular $\|\Seps''\|_{L^\infty} = \|\Weps\|_{L^\infty} = \frac{1}{2\varepsilon}$. Now, assuming particles do not overlap, a simple computation shows for all $i=0,\ldots,N-1$,
\begin{align*}
    & \frac{1}{N}\sum_{k=0}^{N-1}\frac{\Neps(x_{k+1}-x_i)-\Neps(x_k-x_i)}{x_{k+1}-x_k} = \frac{1}{2\varepsilon^2 N}\sum_{k=0}^{N-1}\frac{\left(|x_k-x_i|-|x_{k+1}-x_i|\right)}{x_{k+1}-x_k}\\
    & \ = \frac{1}{2\varepsilon^2 N}\left(2i-N+1\right) =  -\frac{1}{N}\sum_{k=0\,,\,k\neq i}^{N-1}\Neps'(x_i-x_k)-\frac{1}{2\varepsilon^2 N}
\end{align*}
Thus, using difference quotients of $N_\varepsilon$ or using its derivative (avoiding self interactions) in the particle scheme does not change the model in any substantial manner. Therefore, the aforementioned results in \cite{BCDP} imply that the dynamical system
\[\dot{x}_i = \frac{1}{N}\sum_{k=0}^{N-1}\frac{\Neps(x_{k+1}-x_i)-\Neps(x_k-x_i)}{x_{k+1}-x_k}\]
can be formally considered as the gradient flow of the singular, repulsive energy
\[
    \mathcal{N}_\varepsilon[x]=\frac{1}{2 N^2}\sum_{i=0}^{N-1}\sum_{j=0}^{N-1}\Neps(x_i-x_j)-\frac{1}{2\varepsilon^2 N}\sum_{i=0}^{N-1}x_i.
\]
Moreover, using \eqref{eq:splitting}, our particle system \eqref{eq:main_particles} can be formally written as
\[
    \dot{x}_i=\frac{1}{N}\sum_{k=0}^{N-1}\frac{\Neps(x_{k+1}-x_i)-\Neps(x_k-x_i)}{x_{k+1}-x_k}+\frac{1}{N}\sum_{k=0}^{N-1}\frac{\Seps(x_{k+1}-x_i)-\Seps(x_k-x_i)}{x_{k+1}-x_k},
\]
where the first sum can be regarded as the subgradient of $\mathcal{N}_\varepsilon$, and the second sum is the $i$th component of a Lipschitz continuous operator from $\R^{N+1}\rightarrow \R^{N+1}$. The latter can be seen by defining
\[\mathcal{S}_\varepsilon(x,y,z)=\frac{\Seps(x-z)-\Seps(y-z)}{x-y},\]
by observing
\begin{align*}
    & \left|\mathcal{S}_\varepsilon(x,y,z)-\mathcal{S}_\varepsilon(\bar{x},\bar{y},\bar{z})\right|\leq \|\nabla \mathcal{S}_\varepsilon\|_{L^\infty}\|(x,y,z)-(\bar{x},\bar{y},\bar{z})\|\,,
\end{align*}
and by computing for example
\begin{align*}
    & \left|\frac{\partial \mathcal{S}_\varepsilon}{\partial x}\right| = \left|\frac{(x-y)\Seps'(x-z)-\Seps(x-z)+\Seps(y-z)}{(x-y)^2}\right|\\
    & \ = \left|\frac{(x-y)\Seps'(x-z)-\Seps(x-z)+\Seps(x-z)+\Seps'(x-z)(y-x)+\frac{\Seps''(\xi)}{2}(y-x)^2}{(x-y)^2}\right|\leq \frac{1}{4\varepsilon}\,,
\end{align*}
for some suitable intermediate point $\xi\in \R$. Similar estimates for $\frac{\partial \mathcal{S}}{\partial y}$ and $\frac{\partial \mathcal{S}}{\partial z}$ can be obtained easily.

Our goal is to apply Theorem \ref{thm:brezis}. As in Subsection \ref{subsec:gradient_flow_solutions}, to ensure that particle stay in the convex cone $K_N$ for all times we define
the indicator function $\mathcal{I}_{K_N}:\R^{N+1}\rightarrow [0,+\infty]$ as
\[\mathcal{I}_{K_N}(x) = 
\begin{cases}
0 & \hbox{if $x\in K_N$}\\
+\infty & \hbox{if $x\not\in K_N$}\,.
\end{cases}
\]
The above motivates the following definition of solution to \eqref{eq:scheme} with initial conditions in $K_N$.

\begin{defn}\label{def:GF_finite_dimentional}
    Let $\overline{x}\in K_N$. We say that a curve $x:[0,
+\infty)\rightarrow \R^{N+1}$ is a \emph{gradient flow solution} to \eqref{eq:scheme} with initial condition $\overline{x}$ if $x$ is  Lipschitz continuous on $[0,+\infty)$ and satisfies $x(0)=\overline{x}$ and, for all $t>0$,
\begin{equation}\label{eq:GF_finite_dimensions}
    -\dot{x}(t)+\left(\frac{1}{N}\sum_{k=0}^{N-1}\frac{\Seps(x_{k+1}-x_i)-\Seps(x_k-x_i)}{x_{k+1}-x_k}\right)_{i=0}^N\in \partial_{x(t)}\left(\mathcal{N}_\varepsilon+\mathcal{I}_{K_N}\right)\,.
\end{equation}
\end{defn}
We observe that the $S_\varepsilon$ difference quotients in \eqref{eq:GF_finite_dimensions} are well defined even in the case of overlapping particles. Indeed, they can be easily extended by continuity to the whole $K_N$ due to the regularity of $S_\varepsilon$. Hence, \eqref{eq:GF_finite_dimensions} is well defined as a dynamical system in $K_N$. We can therefore state the following theorem.

\begin{thm}\label{thm:finite_dimensional_GF}
Let $\overline{x}\in K_N$. Then, there exists one and only one gradient flow solution to \eqref{eq:scheme} with initial condition $\overline{x}$ in the sense of Definition \ref{def:GF_finite_dimentional}. Moreover, for all $t>0$ we have $x(t)\in O_N$ and $x$ satisfies \eqref{eq:scheme} in a classical sense. Finally, given two gradient flow solutions $x_1(t)$ and $x_2(t)$ to \eqref{eq:scheme}, there exists a positive constant $C_\varepsilon>0$ depending only on $\varepsilon$ such that
\begin{equation}\label{eq:GF_stability}
    \|x_1(t)-x_2(t)\|_N\leq e^{C_\varepsilon t}\|x_1(0)-x_2(0)\|_N
\end{equation}
for all $t\geq 0$, where $\|\cdot\|_N$ is the scaled Euclidean norm on $\R^{N+1}$ introduced in \eqref{eq:scaled-norm}.
\end{thm}

\begin{proof}
Let us define the functional
\[\mathcal{R}_\varepsilon[x]= \frac{1}{4\varepsilon}\left[1-\frac{2}{\varepsilon N}\sum_{i=0}^{N-1}\left(\frac{2i-N}{N}\right) x_i\right],\]    
and observe that $\mathcal{R}_\varepsilon$ is linear and hence convex on $\R^{N+1}$. Therefore, Theorem \ref{thm:brezis} provides the existence and uniqueness of a Lipschitz curve $x:[0,
+\infty)\rightarrow \R^{N+1}$ such that $x(0)=\overline{x}$ and, for $t>0$,
\begin{equation}\label{eq:GF_finite_dimensions_auxiliary}
    -\dot{x}(t)+\left(\frac{1}{N}\sum_{k=0}^{N-1}\frac{\Seps(x_{k+1}-x_i)-\Seps(x_k-x_i)}{x_{k+1}-x_k}\right)_{i=0}^N\in \partial_{x(t)}\left(\mathcal{R}_\varepsilon+\mathcal{I}_{K_N}\right)\,.
\end{equation}
A computation parallel to the one in Lemma \ref{lem:N_and_R_coincide} shows, for all $x\in \R^{N+1}$,
\[\mathcal{R}_\varepsilon[x]+\mathcal{I}_{K_N}[x]=\mathcal{N}_\varepsilon[x]+\mathcal{I}_{K_N}[x].\]
Hence, the existence and uniqueness for \eqref{eq:GF_finite_dimensions} follows. The stability estimate \eqref{eq:GF_stability} is a consequence of \eqref{eq:brezis_stability} using our previous estimates on the Lipschitz map
\[\R^{N+1}\ni x \mapsto \left(\frac{1}{N}\sum_{k=0}^{N-1}\frac{\Seps(x_{k+1}-x_i)-\Seps(x_k-x_i)}{x_{k+1}-x_k}\right)_{i=0}^N\,.\]
Finally, we show that the solution $x(t)$  belongs to $O_N$ for positive times. This would imply automatically that all particles occupy distinct positions for positive times and therefore \eqref{eq:scheme} is satisfied. 
    In case $\overline{x}\in O_N$ then the unique solution $x(t)$ coincides with the one obtained in Subsection \ref{subsec:scheme_strictly_ordered}, for which Lemma \ref{lem:particles_do_not_collide} holds. In particular, estimate \eqref{eq:estimate_d_i} implies that $d_i(t)>0$ for all $i=0,\ldots,N-1$ and all $t>0$. Assume now $\overline{x}\in \partial K_N$, that is there exists $i\in \{0,\ldots,N-1\}$ such that $\overline{x}_i=\overline{x}_{i+1}$. We claim that no particles overlap at positive times $t>0$. We prove the assertion by contradiction. Assume there exists a positive time $t^*>0$ and an index $I\in \{0,\ldots,N-1\}$ such that $x_I(t^*)=x_{I+1}(t^*)$. We consider a new initial condition $\overline{y}=(\overline{y}_0,\ldots,\overline{y}_{N})$ such that $\overline{y}\in O_N$ and with 
    \[
        \|\overline{x}-\overline{y}\|_N<\sigma := \frac{\varepsilon}{2N \sqrt{N+1}}\left(1-e^{-\frac{t^*}{2 \varepsilon^3}}\right)e^{-C_\varepsilon t^*},
    \]
    i.e., $\sigma$ depends on $\varepsilon$, $t^*$, and the constant $C_\epsilon>0$ from the stability estimate.
    
    Since $\overline{y}\in O_N$, \eqref{eq:estimate_d_i} applies to the solution $y(t)$ to \eqref{eq:GF_finite_dimensions} with $y(0)=\overline{y}$. By setting $\widetilde{d}_i=y_{i+1}-y_i$, we therefore get for all $i\in \{0,\ldots,N-1\}$
    \[\widetilde{d}_i(t^*)\geq \frac{2\varepsilon}{N}\left[1-e^{-\frac{t^*}{2\varepsilon^3}}\right]\,.\]
    The stability estimate \eqref{eq:GF_stability} gives
    \[\|x(t^*)-y(t^*)\|_N\leq e^{C_\varepsilon t^*}\|\overline{x}-\overline{y}\|_N < \sigma e^{C_\varepsilon t^*}\,.\]
    Hence,
    \begin{align*}
        & \frac{2\varepsilon}{N}\left[1-e^{-\frac{t^*}{2\varepsilon^3}}\right]\leq |y_{I+1}(t^*)-y_I(t^*)|\\& \ \leq |y_{I+1}(t^*)-x_{I+1}(t^*)| + |y_I(t^*)-x_I(t^*)|+|x_I(t^*)-x_{I+1}(t^*)|\\
        & \ = |y_{I+1}(t^*)-x_{I+1}(t^*)| + |y_I(t^*)-x_I(t^*)| \leq 2 \sqrt{N+1}\sigma e^{C_\varepsilon t^*}\,,
    \end{align*}
    which implies
\[\sigma\geq \frac{\varepsilon}{N \sqrt{N+1}}\left(1-e^{-\frac{t^*}{2 \varepsilon^3}}\right)e^{-C_\varepsilon t^*}\,,\]
a contradiction.
\end{proof}

An important consequence of Theorem \ref{thm:finite_dimensional_GF} is that the unique gradient flow solution satisfies \eqref{eq:scheme_rewritten} for positive times even in the case of initial condition in $K_N$. Clearly, the discrete density $\rho_N$ defined in \eqref{eq:discrete_density} does not make sense at $t=0$ for an initial condition $\overline{x}\in \partial K_N$, that is with overlapping particles. Therefore, the estimates provided in Lemma \ref{lem:estimates_1} lose their meaning in this case. However, estimate \eqref{eq:estimate_d_i} is still meaningful and we can use it to prove a key property of this scheme, which is a discretised version of a \emph{measure-to-$L^p$ smoothing effect}.

\begin{prop}[Smoothing effect]\label{prop:discrete_smoothing}
    Let $x$ be the unique gradient flow solution to \eqref{eq:scheme} with initial condition $\overline{x}$ provided in Theorem \ref{thm:finite_dimensional_GF}. Then, for all $t>0$ and for all $p\in [1,+\infty]$, the discrete density defined in \eqref{eq:discrete_density} satisfies
    \begin{equation}\label{eq:discrete_smoothing}
        \|\rho^N(t)\|_{L^p}\leq\frac{1}{(2\varepsilon)^{\frac{p-1}{p}}}\bigg[ 1-e^{-\frac{t}{2\varepsilon^3}} \bigg]^{-\frac{p-1}{p}}\,.
    \end{equation}
\end{prop}

\begin{proof}
As \eqref{eq:scheme_rewritten} is valid for $t>0$, we still have the validity of \eqref{eq:di_dot} for positive times. Therefore, we can integrate \eqref{eq:di_dot} on $[\delta,t]$ with $0<\delta<t$, and obtain, for all $i=0,\ldots,N-1$, the following estimate for the differences $d_i(t)$ defined in \eqref{eq:differences}: 
\begin{align*}
    d_i (t) \geq \ & e^{-\frac{t}{2\varepsilon^3}} \bigg[ d_i (\delta) + \int_\delta^t \frac{1}{N\varepsilon^2} e^{\frac{s}{2\varepsilon^3}}\,ds \bigg] \\
    = \ & d_i (\delta) e^{-\frac{t}{2\varepsilon^3}} + \frac{2\varepsilon}{N}  \bigg[ e^{-\frac{\delta}{2\varepsilon^3}}-e^{-\frac{t}{2\varepsilon^3}} \bigg]\geq \frac{2\varepsilon}{N}  \bigg[ e^{-\frac{\delta}{2\varepsilon^3}}-e^{-\frac{t}{2\varepsilon^3}} \bigg]\,.
\end{align*}
By sending $\delta\searrow 0$, we obtain
\begin{align*}
    & d_i (t) \geq \frac{2\varepsilon}{N}  \bigg[ 1-e^{-\frac{t}{2\varepsilon^3}} \bigg]\,.
\end{align*}
Hence, the discrete densities $R_i(t)$ defined in \eqref{eq:Ri} satisfy
\[R_i(t)\leq \frac{1}{2\varepsilon}\bigg[ 1-e^{-\frac{t}{2\varepsilon^3}} \bigg]^{-1}\]
and consequently the estimate \eqref{eq:discrete_smoothing} follows in the case $p=+\infty$. The general case $p\in [1,+\infty)$ follows by a standard $L^p$ interpolation recalling that $\rho_N$ has unit $L^1$ norm.
\end{proof}

\subsection{Further estimates on the scheme}\label{subsec:uniform_estimates}

So far, we have proven estimates on the approximated density $\rho^N$ which are uniform with respect to $N$ but depend on $\varepsilon$. While these estimates are very important to analyse the $N\rightarrow+\infty$ limit for fixed $\varepsilon$, they cannot be used in the $\varepsilon\searrow 0$ limit. As it turns out, some of the previous estimates can be improved to be uniform with respect to $\varepsilon$, as we will prove in the present subsection.

Let us consider a general $C^1$-function $\varphi:[0,+\infty)\rightarrow \R$ and set
\[\psi(\rho):=\frac{\varphi(\rho)}{\rho}\,,\]
such that
\[\psi'(\rho)=\frac{\varphi'(\rho)\rho-\varphi(\rho)}{\rho^2}\,.\]

We have the following main lemma.

\begin{lem}[Main uniform estimate]\label{lem:uniform_estimates}
    Assume $\varphi:[0,+\infty)\rightarrow \R$ is $C^1$. Then, the approximated density $\rho^N$ satisfies the estimate
\begin{align}\label{eq:main_uniform_estimate}
    \frac{d}{dt} \int \varphi(\rho^N(y,t))dy = \frac{1}{\varepsilon^2}\int \rho^N(y,t)^2\psi'(\rho^N(y,t))\left[\Weps\ast\rho^N(y,t)-\rho^N(y,t)\right]dy\,.
\end{align}
\end{lem}

\begin{proof}
Let us consider
\begin{align*}
    \frac{d}{dt} &\int \varphi(\rho^N(y,t))dy\\
    &= \frac{d}{dt} \sum_{i=0}^{N-1}\int_{x_i(t)}^{x_{i+1}(t)} \varphi(R_i(t))dy\\
    &= \sum_{i=0}^{N-1} \frac{d}{dt} \left[\varphi(R_i(t)) d_i(t)\right]\\
    &= \sum_{i=0}^{N-1} \left[d_i(t) \varphi'(R_i(t)) \dot R_i(t) + \varphi(R_i(t)) \dot d_i(t)\right]\\
    &= \sum_{i=0}^{N-1}\left[\varphi'(R_i(t))  R_i(t)\int_{x_i(t)}^{x_{i+1}(t)} \Weps''\ast\rho^N(y,t)dy - \varphi(R_i(t))  \int_{x_i(t)}^{x_{i+1}(t)} \Weps''\ast \rho^N(y,t)dy\right],
\end{align*}
having used \eqref{eq:d_i_dot} and \eqref{eq:R_i_dot}.
We can simplify the last expression further to obtain
\begin{align}
    &\sum_{i=0}^{N-1}\left(\int_{x_i(t)}^{x_{i+1}(t)}\rho^N(y,t) \Weps'' \ast \rho^N(y,t)dy\right)\times\left(\varphi'(R_i(t)) - \frac{\varphi(R_i(t))}{R_i(t)}\right)\\
    &= \int \rho^N(y,t) \Weps''\ast\rho^N(y,t) \left( \varphi'(\rho^N(y,t)) - \frac{\varphi(\rho^N(y,t))}{\rho^N(y,t)} \right)dy\,.
\end{align}
Hence, \eqref{eq:morse_elliptic} implies
\[
\frac{d}{dt} \int \varphi(\rho^N(y,t))dy = \frac{1}{\varepsilon^2}\int \rho^N(y,t)^2\psi'(\rho^N(y,t))\left(\Weps\ast\rho^N(y,t)-\rho^N(y,t)\right) dy,
\]
as claimed.
\end{proof}

\begin{cor}\label{cor:rho_Lp}
    For all $N\in \mathbb{N}$ and for all $t\geq 0$ and for all $p\in(1,+\infty]$, the approximated density $\rho^N$ satisfies the estimate
    \begin{equation}\label{eq:rho_Lp}
    \|\rho^N(\cdot,t)\|_{L^p(\R)}\leq \|\rho^N(\cdot,0)\|_{L^p(\R)}.
    \end{equation}
\end{cor}

\begin{proof}
In the case $p\in(1,\infty)$ we consider, with the notation of Lemma \ref{lem:estimates_1},
\[\varphi(\rho)=\rho^p\,,\quad \psi(\rho)=\rho^{p-1},\qquad \psi'(\rho)=(p-1)\rho^{p-2},\]
and use H\"older's inequality and Young's inequality for convolutions to get
\begin{align*}
  \frac{d}{dt} \int_\R (\rho^N(y,t))^p dy &=\frac{p-1}{\varepsilon^2}\int_\R\big[-(\rho^N)(y,t)^{p+1} +(\rho^N(y,t))^p\Weps\ast\rho^N(y,t)\big]dy\\
  & \leq -\frac{p-1}{\varepsilon^2}\int_\R (\rho^N(y,t))^{p+1}dy + \frac{p-1}{\varepsilon^2}\int_\R \Weps dz \int_\R (\rho^N(y,t))^{p+1}dy\\
  &= 0.
\end{align*}
Integrating in time, we get
\begin{equation}\label{eq:rhop}
  \int_\R \rho^N(y,t)^p dy \leq \int_\R \rho^N(y,0)^p dy\,.
\end{equation}
The case $p=+\infty$ follows from the interpolation inequality
\[\|\rho^N(\cdot,0)\|_{L^p(\R)}\leq \|\rho^N(\cdot,0)\|_{L^1(\R)}^{1/p}\|\rho^N(\cdot,0)\|_{L^\infty(\R)}^{(p-1)/p} =\|\rho^N(\cdot,0)\|_{L^\infty(\R)}^{(p-1)/p} \]
and from the inequality
\begin{align*}
\|\rho^N(\cdot,t)\|_{L^\infty(\R)} \leq \limsup_{p\rightarrow+\infty}\|\rho(\cdot,t)\|_{L^p(\R)}
\end{align*}
which follows from
\begin{align*}
    & \|\rho(\cdot,t)\|_{L^p(\R)} \geq \left[\int_{S^N_\delta}\left(\|\rho^N(\cdot,t)\|_{L^\infty(\R)}-\delta\right)^p dx\right]^{1/p} \geq |S_\delta^N|^{1/p}\left(\|\rho^N(\cdot,)\|_{L^\infty(\R)}-\delta\right)
\end{align*}
with $S_\delta^N =\{x\,:\,\, |\rho^N(x,t)|\geq \|\rho^N(\cdot,t)\|_{L^\infty}-\delta\,\,\,\hbox{for a.e. $x\in \R$}\}$, where the $p\rightarrow +\infty$ limit can be taken for fixed $N$ using that $S_\delta^N$ has finite measure (since $\rho^N \in L^1(\R)$) and the assertion follows from the arbitrariness of $\delta>0$.
\end{proof}

\begin{cor}\label{cor:rho_logrho}
    For all $N\in \mathbb{N}$ and for all $t\geq 0$, the approximated density $\rho^N$ satisfies the estimate
    \begin{equation}\label{eq:rho_logrho1}
        \int_\R \rho^N(x,t) \log \rho^N(x,t) dx - \int_0^t \int_\R\rho^N(x,\tau)\Weps''\ast\rho^N(x,\tau)dxd\tau =  \int_\R \rho^N(x,0) \log \rho^N(x,0) dx\,.
    \end{equation}
   Moreover, for all $N\in \mathbb{N}$ and for all $t\geq 0$
     \begin{equation}\label{eq:rho_logrho2}
        \int_\R \rho^N(x,t) \log \rho^N(x,t) dx \leq \int_\R \rho^N(x,0) \log \rho^N(x,0) dx\,.
    \end{equation}
\end{cor}

\begin{proof}
    Consider Lemma \ref{lem:uniform_estimates} with $\varphi(\rho)=\rho \log \rho$. In this case
\[\psi(\rho)=\log\rho\,,\qquad \psi'(\rho)=\frac{1}{\rho}\]
and the estimate \eqref{eq:main_uniform_estimate} reads
\begin{align*}
  & \frac{d}{dt} \int_\R \rho^N \log\rho^N dx =\frac{1}{\varepsilon^2}\int_\R\left[-(\rho^N)^2 +\rho^N\Weps\ast\rho^N\right]dx\\
  & \ =\int_\R \rho^N(x,t)\Weps''\ast\rho^N(x,t)dx.
\end{align*}
Then, estimate \eqref{eq:rho_logrho1} ensues integrating in time. 
Moreover, Young's inequality for convolutions then implies
\begin{align*}
  &  \frac{1}{\varepsilon^2}\int_\R\left[-(\rho^N)^2 +\rho^N\Weps\ast\rho^N\right]dx\leq-\frac{1}{\varepsilon^2}\int_\R(\rho^N)^2 dx + \frac{1}{\varepsilon^2}\int_\R \Weps dx \int_\R(\rho^N)^2 dx \leq 0\,.
\end{align*}
Integrating in time we get \eqref{eq:rho_logrho2}. 
\end{proof}

The estimates proven in Lemma \ref{lem:uniform_estimates} and Corollaries \ref{cor:rho_Lp} and \ref{cor:rho_logrho} show that our particle scheme is well-suited to produce, at the discrete level, estimates for functional of the form $\int\varphi(\rho) dx$ (for a fairly large class of nonlinearities $\varphi$) which one expects to hold at the continuum level too. Such a structure preserving feature of our scheme does not hold as neatly in the case of functionals of the form $\int \rho \phi(x) dx$, which however are extremely useful in the study of the convergence and consistence of the particle scheme, both in the many particle limit for fixed $\varepsilon>0$ and in the joint limit. We prove in the next proposition that 
the evolution of said functionals along the solution to the particle scheme differs from the same computation at the continuum level by an (explicitly computable) error term that vanishes in suitable limiting regimes.

\begin{prop}\label{prop:main_moment_estimate} Let $\rho^N$ be defined as in \eqref{eq:discrete_density} and let $\phi \in C^1 (\R)$ be an arbitrary test function. Then it holds
\[
\frac{d}{dt} \int_\R \phi (x) \rho^N (x,t)\,dx = - \int_\R \rho^N (x,t) \phi' (s) \Weps'\ast \rho^N (x,t)\,dx + C_{\varepsilon,N},
\]
and the error term $C_{\varepsilon,N}$ can be estimated as
\[
\abs{C_{\varepsilon,N}} \leq \frac{\norm{\phi '}_{L^\infty}}{N\varepsilon^2}.
\]
\end{prop}
\begin{proof}
    Let $\phi\in C^1(\R)$ be arbitrary. Recalling \eqref{eq:d_i_dot} and \eqref{eq:R_i_dot} we get 
\begin{align*}
    \frac{d}{dt}\int_\R \phi(x) \rho^N(x, t) dx&=\frac{d}{dt}\sum_{i=0}^{N-1}\int_{x_i(t)}^{x_{i+1}(t)}\phi(x)R_i(t)dx\\
    &=\frac{d}{dt} \sum_{i=0}^{N-1} R_i(t) [\Phi(x_{i+1})- \Phi(x_i)],
\end{align*}
having set $\Phi$ such that $\Phi'(x) = \phi(x)$. Using the product rule, we find
\begin{align}
    \label{eq:time-evo-weak-form}
    \frac{d}{dt}\int_\R \phi(x) \rho^N(x, t) dx &= \mathcal I_1 + \mathcal I_2,
\end{align}
with
\begin{align}
    \label{eq:WeakForm_I1}
    \mathcal I_1 = \sum_{i=0}^{N-1} \dot R_i(t) [\Phi(x_{i+1}(t))- \Phi(x_i(t))],
\end{align}
and
\begin{align}
    \label{eq:WeakForm_I2}
    \mathcal I_2 = \sum_{i=0}^{N-1} R_i(t) \left[
    \phi(x_{i+1}(t))\dot x_{i+1}(t) - \phi(x_i(t)) \dot x_i(t)
    \right].
\end{align}
We continue by treating both terms individually. We omit the time dependence for notation convenience.

\noindent 
\underline{Treatment of $\mathcal I_1$:}  First, let us substitute \eqref{eq:R_i_dot} into \eqref{eq:WeakForm_I1}, to obtain
\begin{align*}
    \mathcal I_1 &= \sum_{i=0}^{N-1} \left\{R_i \fint_{x_i}^{x_{i+1}} \Weps''\ast\rho^N(y,t) dy\right\} [\Phi(x_{i+1})- \Phi(x_i)]\\
    &=\sum_{i=0}^{N-1} R_i \int_{x_i}^{x_{i+1}} \Weps''\ast\rho^N(y,t) \frac{\Phi(x_{i+1})- \Phi(x_i)}{x_{i+1} - x_i}dy \\
    &=\sum_{i=0}^{N-1} \int_{x_i}^{x_{i+1}} \rho^N(y,t)  \Weps''\ast\rho^N(y,t) \frac{\Phi(x_{i+1})- \Phi(x_i)}{x_{i+1} - x_i}dy.
\end{align*}

\noindent 
\underline{Treatment of $\mathcal I_2$:}  Concerning the second contribution, we substitute the particle scheme \eqref{eq:scheme_rewritten} into \eqref{eq:WeakForm_I2} to obtain
\begin{align*}
    \mathcal I_2 &= - \sum_{i=0}^{N-1} R_i(t) \left[\phi(x_{i+1} ) \Weps'\ast\rho^N(x_{i+1}) - \phi(x_i)\Weps'\ast\rho^N(x_i)\right]\\
    &= - \sum_{i=0}^{N-1} R_i(t) \int_{x_i}^{x_{i+1}} \partial_y [\phi(y) \Weps'\ast\rho^N(y,t)] dy,
\end{align*}
having used the fundamental theorem, again. Applying the chain rule, we get
\begin{align*}
    \mathcal I_2 &= - \sum_{i=0}^{N-1} R_i(t) \int_{x_i}^{x_{i+1}} \left[\partial_y \phi(y) \Weps'\ast\rho^N(y,t) + \phi(y) \Weps'' \ast \rho^N(y,t)\right] dy\\
    &= - \sum_{i=0}^{N-1} \int_{x_i}^{x_{i+1}} \rho^N(y,t) \partial_y \phi(y) \Weps'\ast\rho^N(y,t) dy - \sum_{i=0}^{N-1}\int_{x_i}^{x_{i+1}} \phi(y) \rho^N(y,t) \Weps''\ast\rho^N(y,t) dy \\
    &= - \int_\R  \rho^N(y,t) \partial_y \phi(y) \Weps'\ast \rho^N(y,t) dy -  \sum_{i=0}^{N-1} \int_{x_i}^{x_{i+1}} \phi(y) \rho^N(y,t)\Weps''\ast \rho^N(y,t)dy.
\end{align*}
We immediately identify the first term as appearing in the weak formulation in \eqref{eq:weak_nonlocal}.

\noindent
\underline{Combination of $\mathcal I_1$ and $\mathcal I_2$:} Using the expressions we obtained for $\mathcal I_1$ and $\mathcal I_2$, \eqref{eq:time-evo-weak-form} becomes
\begin{align}
\label{eq:weak_error}
    \frac{d}{dt} \int \phi(x) \rho^N(x,t)dx = - \int_\R \rho^N(x,t) \phi'(x) \Weps'\ast \rho^N(x,t) dx + C_{\varepsilon, N},
\end{align} 
where the error term reads
\begin{align}
\label{eq:error1}
    C_{\varepsilon, N} =  \sum_{i=0}^{N-1}\int_{x_i}^{x_{i+1}} \rho^N(s)\Weps''\ast \rho^N(y,t) \left[\frac{\Phi(x_{i+1}) - \Phi(x_i)}{x_{i+1} - x_i} - \Phi'(y)\right] dy,
\end{align}
having used $\Phi'(y) = \phi(y)$, by definition. Next, for $s\in(x_i, x_{i+1})$, a short Taylor expansion reveals
\begin{align*}
    \left[\frac{\Phi(x_{i+1}) - \Phi(x_i)}{x_{i+1} - x_i} - \Phi'(s)\right]  
    &= \left|\frac{\Phi(s + (x_{i+1}-s)) - \Phi(s + (x_i-s))}{x_{i+1} - x_i} - \Phi'(s)\right|\\  
    &\leq \frac12 \norm{\Phi''}_{\infty} (x_{i+1} - x_{i}).
\end{align*}
Hence, we can bound the error in the following way:
\begin{align*}
    |C_{\varepsilon,N}| 
    &\leq \frac12 \|\phi'\|_{L^\infty}  \sum_{i=0}^{N-1} d_i \int_{x_i}^{x_{i+1}} \rho^N(y) |\Weps''\ast\rho^N(y)|  dy\\
    &\leq \frac12 \|\phi'\|_{L^\infty}  \sum_{i=0}^{N-1} d_i  R_i \int_{x_i}^{x_{i+1}} |\Weps''\ast\rho^N(y)|  dy\\
    &\leq \frac1{2N} \|\phi'\|_{L^\infty}  \sum_{i=0}^{N-1} \int_{x_i}^{x_{i+1}} |\Weps''\ast\rho^N(y)|  dy\\
    &= \frac1{2N} \|\phi'\|_{L^\infty}  \int_\R |\Weps''\ast\rho^N(y)|  dy.
\end{align*}
We note that, due to \eqref{eq:morse_elliptic},
\begin{align*}
    \int_\R |\Weps'' \ast \rho^N(y)| dy \leq \frac{1}{\varepsilon^2} \int_\R \left(\rho^N + \Weps \ast \rho^N(y)\right) dy \leq \frac{2}{\varepsilon^2},
\end{align*}
whence
\begin{equation}\label{eq:error}
    |C_{\varepsilon, N}| \leq \frac{\norm{\phi'}_{L^\infty}}{N\varepsilon^2},
\end{equation}
that proves the statement.
\end{proof}

\section{The nonlocal many-particle limit for fixed $\varepsilon>0$}\label{sec:many_particle_limit}

This section is dedicated to proving convergence of the particle scheme \eqref{eq:scheme} as $N\rightarrow +\infty$ and for fixed $\varepsilon>0$. For simplicity, we shall first tackle this task for an initial datum in $L^p$, $p>1$, and then extend the result to initial data in the space of probability measures $\mathcal{P}_2 (\R)$ in the second step. 

\subsection{Definition of the approximation scheme and approximation of the initial data}

We define our approximation scheme in two cases:
\begin{itemize}
    \item [(i)] $\overline{\rho}\in \mpr\cap L^p(\R)$ for $p > 1$;
    \item [(ii)] $\overline{\rho}\in \mathcal{P}_2 (\R)$.
\end{itemize}
The scheme essentially consists of the following three steps:
\begin{itemize}
    \item the approximation of the initial datum $\overline{\rho}$ by a piecewise discrete density featuring $N+1$ jumps at suitable particle positions,
    \item the resolution of the particle scheme via the methods developed in Section \ref{sec:scheme},
    \item the subsequent reconstruction of the piecewise constant density at positive times. 
\end{itemize}

Denoting by $x_0$ the leftmost particle and by $x_N$ the rightmost one, the interval $[x_0,x_N]$ should contain all the mass of the approximated density for all times. Hence, the initial atomisation is relatively easy if $\overline{\rho}$ is compactly supported. If this is not the case, we have to first \emph{compactify} the support of $\overline{\rho}$.

\subsubsection{The case of initial data in $L^p$ for $p > 1$.}
Let us start with case (i). Let $p\in (1,+\infty]$ and let $\overline{\rho}\in \mpr\cap L^p(\R)$, and $N\in \mathbb{N}$ be fixed. Our goal is to define an approximation of $\overline{\rho}$ of the form
\[
\overline{\rho}^N(x)=\sum_{i=0}^{N-1} \overline{R}_i\mathbf{1}_{[\overline{x}_i,\overline{x}_{i+1})}(x)\,,\quad \overline{R}_i=\frac{1}{N\overline{d}_i}\,,\quad \overline{d}_i=\overline{x}_{i+1}-\overline{x}_i,
\]
such that $\overline{\rho}^N$ converges to $\overline{\rho}$ weakly in $L^p(\R)$.
To perform this task, we first define an approximation of $\overline{\rho}$ with compact support and same (unit) mass. For fixed $N$, we set
\begin{align}
    & m_N = \int_{-\infty}^{-N}\overline{\rho}(x) dx\,,\qquad M_N = \int_{N}^{+\infty}\overline{\rho}(x) dx\,.\label{eq:mNMN}
\end{align}
Clearly, both $m_N$ and $M_N$ tend to zero as $N\rightarrow+\infty$.
We then set
\begin{equation}\label{eq:rhoN_tildas}
   \widetilde{\rho}^N(x)=
\begin{dcases}
    \overline{\rho}(x)+\frac{m_N+M_N}{2N}, \quad & \hbox{if $x\in [-N,N]$},\\
    0, & \hbox{otherwise}.
\end{dcases} 
\end{equation}
Notice that $\widetilde{\rho}^N$ has unit mass and is supported on $[-N,N]$. We observe that if $\overline{\rho}$ is compactly supported then $\overline{\rho}$ and $\widetilde{\rho}^N$ coincide for large enough $N$, as in this case $m_N=M_N=0$. We now set
\begin{align}
    & \overline{x}^N_0=-N\label{eq:particle_initial_0}
\end{align}
and iteratively set, for $i=0,\ldots,N-1$,
\begin{align}
    & \overline{x}^N_{i+1}=\inf\left\{x\in \R\,:\,\, \int_{\overline{x}_i^N}^x \widetilde{\rho}^N(y) dy \geq \frac{1}{N}\right\}\,.\label{eq:particle_initial_i}
\end{align}
We observe that $\widetilde{\rho}^N$ has mass $1/N$ between two consecutive points $\overline{x}_i^N$ and $\overline{x}^N_{i+1}$ and that $\overline{x}^N_N = N$. 
For simplicity of notation, we shall drop the superscript $N$ in $\overline{x}^N_i$. We use the points $\overline{x}_0,\ldots,\overline{x}_N$ as initial condition \eqref{eq:scheme_initial} of the particle scheme \eqref{eq:scheme}, whence we define the piecewise constant $\rho^N(x,t)$ as in \eqref{eq:discrete_density}.
The assumption $\overline{\rho}\in L^p$ implies the particles are initially strictly ordered. Lemma \ref{lem:particles_do_not_collide} then implies they remain strictly ordered for all times and the solution $x_0(t), \ldots, x_N(t)$ to \eqref{eq:scheme}-\eqref{eq:scheme_initial} exists globally in time. Our goal is to prove that $\rho^N$ converges in a suitable sense to a weak solution of \eqref{eq:main_intro} with initial condition $\overline{\rho}$, in the sense of Definition \ref{def:weak_nonlocal}. Before stating said result in detail, we first analyse the behavior of $\rho^N$ at time zero. 
We recall
\begin{equation}\label{eq:rhoN_initial}
    \overline{\rho}^N(x)=\rho^N(x,0) =\sum_{i=0}^{N-1}\overline{R}_i\mathbf{1}_{[\overline{x}_i,\overline{x}_{i+1})}(x)\,,\qquad \overline{R}_i=\frac{1}{N(\overline{x}_{i+1} -\overline{x}_i)},
\end{equation}
i.e., $\overline{\rho}^N$ is a piecewise constant approximation of $\tilde \rho^N$. We notice that on each interval $[\overline{x}_i,\overline{x}_{i+1}]$ we have
\begin{equation}\label{eq:initial_average}
    \overline{R}_i = \fint_{\overline{x}_i}^{\overline{x}_{i+1}} \widetilde{\rho}^N(y) dy\,.
\end{equation}

The next lemma is essential to control the $L^p$ norms of $\overline{\rho}^N$ uniformly in $N$ in terms of the $L^p$ norm of $\overline{\rho}$.

\begin{lem}[Uniform control of the $L^p$ norms of the initial data]\label{lem:initial_control}
    Let $\overline{\rho}\in L^p(\R)$ for some $p\in [1,+\infty]$. Then there exists a constant $C \geq 0$ depending only on $\|\overline{\rho}\|_{L^p(\R)}$ such that 
    \[\|\overline{\rho}^N\|_{L^p(\R)}\leq C, \qquad \hbox{for all $N\in \mathbb{N}$}\,.\]
\end{lem}

\begin{proof}
    The statement is trivial for $p=1$. Let us now consider the case $p\in (1,+\infty)$. We may write
\begin{align*}
    \norm{\overline{\rho}^N}_{L^p(\R)}^p 
    &= \int_\R (\overline{\rho}^N (x) ) ^p dx = \sum_{i=0}^{N-1}\int_{\overline{x}_i}^{\overline{x}_{i+1}} {\overline{R}_i }^p dx = \sum_{i=0}^{N-1}\int_{\overline{x}_i}^{\overline{x}_{i+1}}\left(\int_{\overline{x}_i}^{\overline{x}_{i+1}}\widetilde{\rho}^N(y) \frac{dy}{\overline{x}_{i+1}-\overline{x}_i}\right)^p dx.
\end{align*}
Thus, we can estimate the $L^p$ norm using Jensen's inequality, to get
\begin{align*}
    \norm{\overline{\rho}^N}_{L^p(\R)}^p 
    &\leq \sum_{i=0}^{N-1} \int_{\overline{x}_i}^{\overline{x}_{i+1}} \int_{\overline{x}_i}^{\overline{x}_{i+1}}(\widetilde{\rho}^N(y))^p \frac{dy}{\overline{x}_{i+1}-\overline{x}_i} dx\\ 
    &= \sum_{i=0}^{N-1} \int_{\overline{x}_i}^{\overline{x}_{i+1}} (\widetilde{\rho}^N(y))^p dy.
\end{align*}
Finally, using the definition of $\tilde \rho^N$, we get
\begin{align*}
    \norm{\overline{\rho}^N}_{L^p(\R)}^p 
    & = \int_\R (\widetilde{\rho}^N(y))^p dy = \int_{-N}^N \left(\overline{\rho}(x) + \frac{m_N+M_N}{2N}\right)^p dx \\
    &\leq 2^p \int_{-N}^N\left[\overline{\rho}(x)^p +\left(\frac{m_N+M_N}{2N}\right)^p\right] dx\\
    & =  2^p\left(\|\overline{\rho}\|_{L^p(\R)}^p + (m_N+M_N)^p(2N)^{1-p}\right)\,.
\end{align*}
The second term is uniformly bounded with respect to $N$ since $\overline{\rho}\in L^p(\R)$. The case $p=+\infty$ follows from \eqref{eq:initial_average} and from the fact that the $L^\infty$ norm of $\widetilde{\rho}^N$ is uniformly bounded in $N$ as a direct consequence of the definition \eqref{eq:rhoN_tildas}.
\end{proof}

We now collect in the next lemma some of the convergence properties of $\overline{\rho}^N$ towards $\overline{\rho}$.
\begin{lem}\label{lem:initial_convergence}
    Let $\overline{\rho}\in\mpr\cap L^p(\R)$ for some $p\in (1,+\infty]$. Then, $\overline{\rho}^N\rightarrow \overline{\rho}$ in the sense of distributions. Moreover, if $p\in (1,+\infty)$ then $\overline{\rho}^N\rightharpoonup \overline{\rho}$ weakly in $L^p(\R)$, otherwise, if $p=+\infty$ the convergence holds in the $L^\infty(\R)$ weak-$*$ sense.
\end{lem}

\begin{proof}
Let $\phi \in C^\infty_c(\R)$ be an arbitrary test function. Let $L>0$ be such that the support of $\phi$ is contained in $[-L,L]$. Let $\overline{x}_I$ be the rightmost particle with position strictly less than $-L$ and $\overline{x}_J$ be the leftmost particle with position strictly larger than $L$. We compute
\begin{align*}
    & \int_\R \left(\overline{\rho}^N(x)-\widetilde{\rho}^N(x)\right) \phi(x) dx=\int_{-L}^L \left(\overline{\rho}^N(x)-\widetilde{\rho}^N(x)\right) \phi(x) dx  \\
    & \qquad = \sum_{i=I}^{J-1}\int_{\overline{x}_i}^{\overline{x}_{i+1}}\left(\overline{\rho}^N(x)-\widetilde{\rho}^N(x)\right) \phi(x) dx \\
    & \qquad = \sum_{i=I}^{J-1}\int_{\overline{x}_i}^{\overline{x}_{i+1}}\phi(x)\left(\fint_{\overline{x}_i}^{\overline{x}_{i+1}}\widetilde{\rho}^N(y)dy-\widetilde{\rho}^N(x)\right) dx \\
    & \qquad =\sum_{i=I}^{J-1}\frac{1}{\overline{x}_{i+1}-\overline{x}_i}\int_{\overline{x}_i}^{\overline{x}_{i+1}}\int_{\overline{x}_i}^{\overline{x}_{i+1}}\phi(x)\left(\widetilde{\rho}^N(y)-\widetilde{\rho}^N(x)\right) dy dx\\
    & \qquad = \sum_{i=I}^{J-1}\frac{1}{\overline{x}_{i+1}-\overline{x}_i}\int_{\overline{x}_i}^{\overline{x}_{i+1}}\int_{\overline{x}_i}^{\overline{x}_{i+1}}\left(\phi(y)-\phi(x)\right)\widetilde{\rho}^N(x) dy dx = A_N + B_N,
\end{align*}
with
\begin{align*}
    & A_N=\frac{1}{\overline{x}_{I+1}-\overline{x}_I}\int_{\overline{x}_I}^{\overline{x}_{I+1}}\int_{\overline{x}_I}^{\overline{x}_{I+1}}\left(\phi(y)-\phi(x)\right)\widetilde{\rho}^N(x) dy dx \\
    & \qquad\quad  +
    \frac{1}{\overline{x}_{J}-\overline{x}_{J-1}}\int_{\overline{x}_{J-1}}^{\overline{x}_{J}}\int_{\overline{x}_{J-1}}^{\overline{x}_{J}}\left(\phi(y)-\phi(x)\right)\widetilde{\rho}^N(x) dy dx,
\end{align*}
and
\begin{align*}
    & B_N = \sum_{i=I+1}^{J-2}\frac{1}{\overline{x}_{i+1}-\overline{x}_i}\int_{\overline{x}_i}^{\overline{x}_{i+1}}\int_{\overline{x}_i}^{\overline{x}_{i+1}}\left(\phi(y)-\phi(x)\right)\widetilde{\rho}^N(x) dy dx.
\end{align*}
The term $A_N$ above is the sum of two terms that can be controlled similarly. Hence, we only present the detailed estimate for the first term, namely
\begin{align*}
    & \left|\frac{1}{\overline{x}_{I+1}-\overline{x}_I}\int_{\overline{x}_I}^{\overline{x}_{I+1}}\int_{\overline{x}_I}^{\overline{x}_{I+1}}\left(\phi(y)-\phi(x)\right)\widetilde{\rho}^N(x) dy dx\right|\\
    &\qquad \leq \frac{2\|\phi\|_\infty}{\overline{x}_{I+1} - \overline x_{I}} (\overline x_{I+1} - \overline{x}_{I}) \int_{\overline x_I}^{\overline x_{I+1}} \widetilde \rho^N(x) dx  =\frac{2\|\phi\|_\infty}{N}.
\end{align*}
Hence, we have proven that $A_N\rightarrow 0$ as $N\rightarrow +\infty$. Using the regularity of $\phi$, we now estimate $B_N$:
\begin{align*}
    \left|B_N\right|
    &\leq \|\phi'\|_{L^\infty}\sum_{i=I+1}^{J-2}(\overline{x}_{i+1}-\overline{x}_i)\int_{\overline{x}_i}^{\overline{x}_{i+1}}\widetilde{\rho}^N(x) dx= \frac{\|\phi'\|_{L^\infty}}{N}\sum_{i=I+1}^{J-2}(\overline{x}_{i+1}-\overline{x}_i) \\
&=  \frac{\|\phi'\|_{L^\infty}}{N} (\overline{x}_{J-1}-\overline{x}_{I+1})\leq \frac{2L\|\phi'\|_{L^\infty}}{N}\,.
\end{align*}
The above estimates prove that the difference $\overline{\rho}^N-\widetilde{\rho}^N$ goes to zero in the sense of distributions. Now, we may compute
\begin{align}
    \left|\int_\R (\widetilde{\rho}^N(x)-\overline{\rho}(x))\phi(x) dx\right| &\leq (m_N + M_N)\|\phi\|_{L^\infty} + \int_{-N}^N \left|\widetilde{\rho}^N(x)-\overline{\rho}(x)\right|\phi(x) dx\nonumber\\
    & = 2(m_N + M_N) \|\phi\|_{L^\infty},\label{eq:estimate_rearrange}
\end{align}
and the latter tends to zero as well as $N\rightarrow +\infty$. We have therefore proven that $\overline{\rho}^N$ tends to $\overline{\rho}$ in the sense of distributions as $N\rightarrow +\infty$. In order to prove the weak $L^p$ convergence, we recall from Lemma \ref{lem:initial_control} that the $L^p$ norm of $\overline{\rho}^N$ is uniformly bounded with respect to $N$. In the case $p>1$ we use the reflexivity of $L^p$ to obtain that a subsequence of $\overline{\rho}^N\in L^p(\R)$ converges weakly. The only possible limit is $\overline{\rho}$ because of the distributional convergence proved above. Hence, the whole sequence $\overline{\rho}^N$ converges weakly to $\overline{\rho}$. The case $p=+\infty$ follows similarly from Banach-Alaoglu's Theorem.
\end{proof}

\subsubsection{The case of measure initial data.}
In order to extend the scheme to the general case (ii) we use an alternative approach based on the use of the quantile function, see Subsection \ref{subsec:one_d_wass}. Let $\overline{\rho}\in \mathcal{P}_2 (\R)$ and $N\in \mathbb{N}$. Similar to the previous case (i), we set
\begin{equation}\label{eq:mNMN2}
     m_N = \overline{\rho}((-\infty,-N))\,,\qquad M_N = \overline{\rho}((N,+\infty))
\end{equation}
and define, for all Borel measurable sets $A\subset \R$,
\begin{equation}\label{eq:rhoN_tildas_2}
    \widetilde{\rho}^N(A)=\overline{\rho}(A\cap[-N,N]) + \frac{m_N}{N}\mathcal{L}^1(A\cap[-N,0))+ \frac{M_N}{N}\mathcal{L}^1(A\cap[0,N])\,,
\end{equation}
where $\mathcal{L}^1$ is the one-dimensional Lebesgue measure. Observe that $\widetilde{\rho}^N$ is a probability measure with compact support. Now, given the isometry $\mathcal{T}$ defined in Subsection \ref{subsec:one_d_wass}, we set
\begin{equation}\label{eq:X_tilda_N}
    \widetilde{X}^N \coloneqq \mathcal{T}(\widetilde{\rho}^N)\in \mathcal{K}\subset L^2((0,1))\,.
\end{equation}
We recall that $\widetilde{X}^N$ is non-decreasing and right-continuous. Moreover, since $\mathrm{supp}(\widetilde{\rho}^N)\subset [-N,N]$, we have the estimate
\[-N\leq \widetilde{X}^N(z) \leq N\,,\qquad \hbox{for all $z\in [0,1]$}\,.\]
In order to atomise $\widetilde{\rho}^N$, we split the interval $[0,1]$ into $N$ intervals of size $1/N$ by setting $z_i=i/N$ for $i=0,\ldots,N$. 

We now define
\begin{equation}\label{eq:atomisation_X_tilde}
    \overline{X}^N(z)=\sum_{i=0}^{N-1}\left[\widetilde{X}^N(z_i) + N\left(\widetilde{X}^N(z_{i+1})-\widetilde{X}^N(z_i)\right) (z-z_i)\right]\mathbf{1}_{[z_i,z_{i+1})}(z)\,.
\end{equation}
In a nutshell, $\overline{X}^N$ is a piecewise linear interpolation of the points $(z_i,\widetilde{X}^N(z_i))\in [0,1]\times [-N, N]$. Note that $\overline{X}^N$ is constant on the intervals $[z_i,z_{i+1}]$ where $\widetilde{X}^N(z_i)=\widetilde{X}^N(z_{i+1})$, which happens when the measure $\widetilde{\rho}^N$ has a mass larger or equal to $1/N$ concentrated on the point $\widetilde{X}^N(z_i)$. This situation occurs for $N$ large enough in case of a nonzero mass concentration for $\overline{\rho}$ on the same point. We finally set
\begin{equation}\label{eq:initial_discrete_measure}
  \overline{\rho}^N=\mathcal{T}^{-1}(\overline{X}^N)\,,
\end{equation}
and
\begin{equation}\label{eq:initial_particles_measure}
 \overline{x}_i=\widetilde{X}^N(z_i)\,,\qquad i=0,\ldots,N\,.   
\end{equation}
We then define the particle scheme as in Subsection \ref{subsec:scheme_overlapping}. Note that particles may overlap initially in this case. More precisely, given
\[\mathcal{S}^N=\left\{i\in\{0,\ldots,N-1\}\,:\,\, \overline{x}_{i+1}=\overline{x}_i\right\},\]
we have
\begin{equation}\label{eq:rhoN_initial_measure}
    \overline{\rho}^N=\frac{1}{N}\sum_{i\in \mathcal{S}^N}\delta_{\overline{x}_i} + \sum_{i\in \{0,\ldots,N-1\}\setminus\mathcal{S}^N}\overline{R}_i\mathbf{1}_{[\overline{x}_i,\overline{x}_{i+1})},
\end{equation}
with the usual notation
\[\overline{R}_i=\frac{1}{N(\overline{x}_{i+1}-\overline{x}_i)}\,.\]
We therefore consider the unique gradient flow solution $x_0(t),\ldots,x_N(t)$ to \eqref{eq:scheme} provided by Theorem \ref{thm:finite_dimensional_GF}. The smoothing effect result in Proposition \ref{prop:discrete_smoothing} guarantees that particles detach immediately for positive times and stay like that globally in time. Hence, we may define the piecewise constant density $\rho^N(x,t)$ as in \eqref{eq:discrete_density}.

\begin{lem}\label{lem:initial_convergence_measure}
    Let $\overline{\rho} \in \mathcal{P}_2 (\R)$. Then, the sequence $\overline{\rho}^N$ converges to $\overline{\rho}$ in the weak measure sense.
\end{lem}
\begin{proof}\mbox{}\\
\noindent
\underline{Step 1 -- Tightness of $\overline{\rho}^N$:}
 The sequence $\overline{\rho}^N$ is \emph{tight} due to Prokhorov's Theorem (see for example \cite{AGS}). To see this, let $\varepsilon>0$. Since $\overline{\rho}$ is a probability measure, there exists a bounded set $B\subset \R$ with $\overline{\rho}(\R\setminus B)<\varepsilon/3$. Then, for all $N\in \mathbb{N}$ we have $\widetilde{\rho}^N(B)\geq \overline{\rho}(B)-m_N-M_N$. Consequently, 
 \[\widetilde{\rho}^N(\R\setminus B)=1-\widetilde{\rho}^N(B)\leq 1-\overline{\rho}(B)+m_N+M_N=\overline{\rho}(\R\setminus B)+m_N+M_N<\varepsilon/2+m_N+M_N\,,\]
and since $m_N+M_N\rightarrow 0$, there exists $N_\varepsilon\in \mathbb{N}$ such that
\[
    \widetilde{\rho}^N(\R\setminus B)<\frac{2}{3}\varepsilon,
\]
for all $N\geq N_\varepsilon$. Now, by definition of $\overline{\rho}^N$, for any Borel set $A\subset \R$ we have
\[\overline{\rho}^N(A)=\mathcal{L}^1((\overline{X}^N)^{-1}(A))\,,\]
and since
\[
    0\leq \left|\overline{X}^N(z)-\widetilde{X}^N(z)\right|\leq \widetilde{X}^N(z_{i+1})-\widetilde{X}^N(z_i),
\]
with $z\in [z_i,z_{i+1})$, recalling that $\mathcal{L}^1([z_{i},z_{i+1}))=1/N$, we get that $\overline{\rho}^N(A)$ and $\widetilde{\rho}^N(A)$ differ at most by $1/N$. Hence, by possibly renaming $N_\varepsilon$ and for $N\geq N_\varepsilon$ we get
\[
    \overline{\rho}^N(\R\setminus B)<\widetilde{\rho}^N(\R\setminus B)+\varepsilon/3<\varepsilon,
\]
for $N\geq N_\varepsilon$. We notice that $N_\varepsilon$ does not depend on $B$, which proves the tightness of $\overline{\rho}^N$.

Then, up to a (non-relabelled) subsequence, $\overline{\rho}^N$ converges in the weak sense of measures to some probability measure $\rho$. It remains to identify the limit as the original initial data.
In order to show that $\rho=\overline{\rho}$, we switch to the associated pseudoinverse variables and
\begin{itemize}
    \item show that $\widetilde X^N$ has a pointwise limit, $X$
    \item the limit coincides with the limit of $\overline{X}^N$
    \item and $X = \overline{X} := \mathcal{T}(\overline{\rho})$, i.e., the limit is identical to the pseudoinverse of the initial data, $\overline \rho$.
\end{itemize}

\noindent
\underline{Step 2 -- Convergence of $\widetilde X^N$:} 
Let $K=[a,b]\subset(0,1)$ be a compact interval. First, we observe that $\widetilde X^N$ is uniformly bounded on $K$ as a consequence of the moment control of the associated measure $\widetilde \rho^N$. Indeed, let us set $x^N_a=\widetilde{X}^N(a)$ and observe that, by definition,
\[
\widetilde{\rho}^N((-\infty,x^N_a))=a>0.
\]
Hence, if $x_a^N$ were to diverge to $-\infty$ (up to a subsequence), we would have a contradiction. Similarly, we can prove that $x^N_b=\widetilde{X}^N(b)$ cannot diverge to $+\infty$ up to a subsequence. Next, due to our construction in \eqref{eq:rhoN_tildas_2}, we observe
\begin{subequations}
\label{eq:monotonicity-tildeXN}
\begin{align}
    & \widetilde{X}^N(z) \geq  \widetilde{X}^{N+1}(z) \geq \overline{X}(z), \qquad \hbox{if $0\leq z<1/2$},
\end{align}
as well as
\begin{align}
    & \widetilde{X}^N(z)\leq \widetilde{X}^{N+1}(z) \leq \overline{X}(z), \qquad \hbox{if $1/2< z\leq 1$},
\end{align}
\end{subequations}
and $y_0 := \overline{X}(1/2)=\widetilde{X}^N(1/2)$ for all $N\in \mathbb{N}$, as \eqref{eq:rhoN_tildas_2} preserves the total masses on $(-\infty,y_0)$ and $[y_0,+\infty)$, and the fact that 
$$
    \partial_z \widetilde{X}^N(z)\leq \partial_z\overline{X}(z)\qquad \text{on } z\in [0,1/2),
$$
as well as 
$$
    \partial_z \widetilde{X}^N(z)\geq \partial_z\overline{X}(z), \qquad \text{on } z\in [1/2,1].
$$ 
Due to the monotonicity properties of the sequence $(\widetilde X^N)_N$, \eqref{eq:monotonicity-tildeXN} , we can set $X(z) := \lim_{N\to \infty} \widetilde X^N(z)$, for each $z\in [0,1]$.

\noindent
\underline{Step 3 -- Limits of $\widetilde X^N$ and $\overline X^N$ coincide:} 
Now, let $L>0$ be large enough to have  $[\widetilde{X}^N(a),\widetilde{X}^N(b)]\subset [-L,L]$ for all $N\in \mathbb{N}$. Let $I\in \{0,\ldots,N\}$ be the smalles index such that $\widetilde{X}^N(z_I)\geq -L$ and let $J\in \{0,\ldots,N\}$ be the largest index such that $\widetilde{X}^N(z_J)\leq L$. For all $\delta>0$ small and fixed and for $N$ large enough we have
    \begin{align*}
        & \int_{a+\delta}^{b-\delta} |\widetilde{X}^N(z)-\overline{X}^N(z)| dz\\
        &\qquad \leq \sum_{i=I}^{J-1}\int_{z_i}^{z_{i+1}} |\widetilde{X}^N(z)-\overline{X}^N(z)| dz\\
        &\qquad  \leq \sum_{i=I}^{J-1}\int_{z_i}^{z_{i+1}} |\widetilde{X}^N(z_{i+1})-\widetilde{X}^N(z_i)| dz,
    \end{align*}
where the last step is due to the fact that $\overline{X}^N$ is the piecewise linear interpolation of $\widetilde{X}^N$. Hence,
\begin{align*}
     & \int_{a+\delta}^{b-\delta} |\widetilde{X}^N(z)-\overline{X}^N(z)| dz\leq \sum_{i=I}^{J-1}\frac{\widetilde{X}^N(z_{i+1})-\widetilde{X}^N(z_i)}{N}=\frac{\widetilde{X}^N(z_{J})-\widetilde{X}^N(z_I)}{N}\leq \frac{2L}{N}\,.
\end{align*}
Since $\delta>0$ is arbitrary, we have therefore proven that $\widetilde{X}^N-\overline{X}^N$ tends to zero in $L^1_{\mathrm{loc}}((0,1))$ and hence almost everywhere on $[0,1]$ up to a subsequence and after a diagonal procedure.

\noindent
\underline{Step 4 -- Identification of $X$:} We argue by contradiction and assume the pointwise limit, $X(z)= \lim_{N\to\infty}\widetilde{X}^N(z)$ is different from $\overline X$. Then, by Lebesgue's dominated convergence, we have for any $\phi\in C_b(\R)$ that
\[
    \int_0^1 \phi(\widetilde{X}^N(z)) dz \rightarrow \int_0^1 \phi(X(z))dz,
\]
and from \eqref{eq:change_of_variable} this implies that $\widetilde{\rho}^N$ converges in the weak measure sense to a measure that is not equal to $\overline{\rho}$. However, with an estimate similar to the one in \eqref{eq:estimate_rearrange} one can easily see that this is a contradiction.  Once again taking $\phi\in C_b(\R)$ yields
\[\int_\R \phi(x) d\overline{\rho}^N(x) = \int_0^1 \phi(\overline{X}^N(z))dz \rightarrow \int_0^1 \phi(\overline{X}(z)) dz = \int_\R \phi(x) d\overline{\rho}(x), \]
which, in conjunction with Step 3 above, shows the assertion.

\end{proof}

\subsection{Convergence of the scheme}\label{subsec:convergence1}

In this section we prove the two main results of convergence of the particle scheme as $N\rightarrow +\infty$ and for fixed $\varepsilon>0$. We shall first state and prove our result with an initial condition in $\overline{\rho}\in L^p(\R)\cap \mathcal{P}_2 (\R)$ for $p>1$ and then in the more general case of an initial condition in $\mathcal{P}_2 (\R)$.

\begin{thm}[Many particle limit for fixed $\varepsilon$ with initial data in $L^p$]\label{thm:convergence1}
    Let $p\in (1,+\infty]$ and let $\overline{\rho}\in\mptr\cap L^p(\R)$. Let $\overline{\rho}^N$ be as in \eqref{eq:rhoN_initial}, with $\overline{x}_0,\ldots,\overline{x}_N$ defined in \eqref{eq:particle_initial_0}-\eqref{eq:particle_initial_i} and $\widetilde{\rho}^N$ defined in \eqref{eq:rhoN_tildas}. 
    Let $x_0(t),\ldots,x_N(t)$ be the unique solution to the scheme \eqref{eq:scheme} with initial datum $\overline{x}_0,\ldots,\overline{x}_N$, and let $\rho^N$ be the corresponding discrete piecewise constant density as in \eqref{eq:discrete_density}. 
    Then, for any fixed $T\geq 0$, $\rho^N$ converges weakly in $L^p([0,T]\times \R)$ for $p\in[-1,\infty)$, (resp. weakly-$^*$ in $L^\infty([0,T]\times \R)$, if $p=\infty$), to the unique weak measure solution $\rho$ to \eqref{eq:main_intro} in the sense of Definition \ref{def:weak_nonlocal}.
\end{thm}

\begin{proof}
\noindent
By using both the uniform  $L^\infty(0,T; L^p(\R))$-bounds on $\rho^N$ provided by Lemma \ref{lem:estimates_1} and the uniform estimate on $\overline{\rho}^N$ stated in Lemma \ref{lem:initial_control}, we may extract a non-relabelled subsequence converging weakly (weakly-$^*$ in the $p=+\infty$ case) to a limit $\rho\in L^p(\R\times [0,T])$ as $N\to \infty$. Concurrently, we have
\begin{align}
    \label{eq:pointwise-conve-convo}
    \lim_{N\to \infty} \abs{\Weps' \ast \rho^N(x,t) - \Weps' \ast \rho(x,t)} \leq \left| \lim_{N\to\infty} \int_\R \Weps'(x-y) (\rho^N(y,t) - \rho(y,t)) dy \right|= 0,
\end{align}
since $\Weps' \in L^q(\R)$, where $q\in [1,+\infty)$ denotes the  Hölder conjugate associated to $p\in (1,+\infty]$. Young's inequality for convolution provides $\norm{\Weps' \ast \rho^N}_{L^\infty} \leq \norm{\Weps'}_{L^q}$, which, in conjunction with the pointwise convergence in \eqref{eq:pointwise-conve-convo}, yields 
\begin{align}
    \label{eq:strong-conve-convo}
    \Weps'\ast \rho^N \rightarrow \Weps' \ast \rho,
\end{align}
strongly in $L^q_{\mathrm{loc}}(\R \times (0,T))$. We recall that Lemma \ref{lem:initial_convergence} implies $\rho^N(\cdot,0) \rightharpoonup \overline{\rho}$, weakly in $L^p(\R)$.
By Proposition \ref{prop:main_moment_estimate} we know that
\begin{align}
    \frac{d}{dt} \int_\R \phi(x) \rho^N(x, t) dx = - \int_\R \rho^N(x, t) \Weps' \ast \rho^N(x, t) \phi'(x)  dx + \mathcal O(N^{-1} \varepsilon^{-2})\,.
\end{align}
Multiplying by another test function, $\chi \in C_c^\infty([0,\infty))$ and integrating over the positive real half line, having picked $T>0$ such that $\mathrm{supp}(\chi)\subset [0,T]$, we get
\begin{align*}
    \int_0^{T} \chi(t) \frac{d}{dt} \int_\R \phi(x) \rho^N(x, t) dx dt  &= - \int_0^{T}\int_\R \rho^N(x, t)  \Weps' \ast \rho^N(x, t) \phi'(x) \chi(t) dx dt + \mathcal O(T N^{-1} \varepsilon^{-2}).
\end{align*}
Integrating by parts in time finally gives us
\begin{align}
    \label{eq:almost-weak-form}
    \begin{split}
    \int_0^T\int_\R \big[\partial_t \varphi(x,t)  -   \Weps' \ast \rho^N(x, t) \partial_x \varphi(x,t) \big] \rho^N(x,t)  dx dt &=  - \int_\R \varphi(x,0) \rho^N(x,0) dx\\
    &\qquad + \mathcal O(T N^{-1} \varepsilon^{-2}).
    \end{split}
\end{align}
where we set $\varphi(x,t) := \chi(t) \phi(x)$. By using a density argument, we have that the above holds for any test function $\varphi = \varphi (x,t)$, with $\varphi\in C^1_c(\R\times [0,T))$, see \cite[Thm.\ 4.3.1]{friedlander}.
Using the (up to a subsequence) weak convergence of $\rho^N$ to $\rho$, the strong convergence of the convolution term in \eqref{eq:strong-conve-convo}, and the weak convergence of $\rho^N(\cdot,0)$ proven in Lemma \ref{lem:initial_convergence}, we obtain
\begin{align}
    \label{eq:weak-form}
    \int_0^T\int_\R \big[\partial_t \varphi(x,t)  -   \Weps' \ast \rho(x, t) \partial_x \varphi(x,t) \big] \rho(x,t)  dx dt = - \int_\R \varphi(x,0) \bar \rho(x) dx,
\end{align}
which proves that $\rho$ is a weak solution in the sense of Definition \ref{def:weak_nonlocal}.

To prove the uniqueness statement, we observe that Proposition \ref{prop:discrete_smoothing} implies that $\rho^N$ is uniformly bounded on $(x,t)\in \R\times [\delta,T]$ for every fixed $\delta>0$. Hence, by possibly extracting another subsequence, $\rho^N$ converges weakly $^*$ in $L^\infty(\R\times [\delta,T])$. The weak lower semi-continuity of the $L^\infty$ norm implies
\[\|\rho\|_{L^\infty(\R\times[\delta,T])}\leq \liminf_{N\rightarrow+\infty}\|\rho^N\|_{L^\infty(\R\times[\delta,T])}\leq C\]
for some suitable $C\geq 0$. Hence, $\rho$ satisfies all the assumptions and the regularity prescribed in Definition \ref{def:weak_nonlocal}.
Statement (v) in Theorem \ref{thm:existence_and_uniqueness_GF} then implies that $\rho$ is a gradient flow solution in the sense of Definition \ref{def:GF}, which is unique due to the same Theorem \ref{thm:existence_and_uniqueness_GF}. Hence, the whole sequence $\rho^N$ converges to $\rho$.
\end{proof}

We now extend the above theorem to provide our main convergence result, which deals with initial data in $\mptr$. 

\begin{thm}[Many particle limit for fixed $\varepsilon$ with measure initial data]\label{thm:convergence_measure}
    Let $\overline{\rho}\in \mptr$. Let $\overline{\rho}^N$ be as in \eqref{eq:initial_discrete_measure}-\eqref{eq:atomisation_X_tilde}-\eqref{eq:X_tilda_N}-\eqref{eq:rhoN_tildas_2}-\eqref{eq:mNMN2}. Let $x_0(t),\ldots,x_N(t)$ be the unique gradient flow solution provided by Theorem \ref{thm:finite_dimensional_GF} with initial data \eqref{eq:initial_particles_measure} and let $\rho^N$ be the corresponding discrete piecewise constant density as in \eqref{eq:discrete_density}. Then, for fixed $T\geq 0$, $\rho^N$ converges weakly $^*$ in $L^\infty_{\mathrm{loc}}((0,T]\,;\,L^\infty(\R))$ to the unique gradient flow solution $\rho$ to \eqref{eq:main_intro} with initial condition $\overline{\rho}$ provided by Theorem \ref{thm:existence_and_uniqueness_GF}.
\end{thm}

\begin{proof}
Theorem \ref{thm:finite_dimensional_GF} implies that particles do not overlap for positive times. Hence, we can apply the same computation performed in the proof of Theorem \ref{thm:convergence1} for times $t\geq \delta$ with $\delta>0$ given and obtain
\begin{align}
    \label{eq:almost-weak-form2}
    \int_\delta^T\int_\R \big[\partial_t \varphi(x,t)  -   \Weps' \ast \rho^N(x, t) \partial_x \varphi(x,t) \big] \rho^N(x,t)  dx dt &=  - \int_\R \varphi(x,\delta) \rho^N(x,\delta) dx\\
    &\qquad + \mathcal O(T N^{-1} \varepsilon^{-2})
\end{align}
for all test functions $\varphi\in C^1_c(\R\times [\delta,T))$. We notice in particular that the estimate of the error term only uses the $L^1$ norm of $\rho^N$, so the above computation still applies to this case and the error is independent of $\delta$. 

We now let $\delta\searrow 0$. We recall
\[\rho^N(x,\delta)=\sum_{i=0}^{N-1}R_i(\delta)\mathbf{1}_{[x_i(\delta),x_{i+1}(\delta))}(x)\]
and observe that for fixed $N$ the above expression tends, as $\delta\searrow 0$ in the $2$-Wassertein sense, to $\overline{\rho}^N$ defined in \eqref{eq:initial_discrete_measure}. To see this, recall that $\overline{\rho}^N$ can be written as in \eqref{eq:rhoN_initial_measure}. The trajectories of the particles $x_i(t)$ for all $i=0,\ldots,N-1$ are Lipschitz continuous on $[0,\delta]$ and tend to $\overline{x}_i$ respectively for all $i=0,\ldots,N-1$ as $\delta\searrow0$. Hence, with the notation in \eqref{eq:rhoN_initial_measure}, we have for all $\varphi\in C(\R)$ with finite second moment that
\begin{align*}
    \int_\R \rho^N(x,\delta) \varphi(x) dx 
    &= \sum_{i=0}^{N-1}\frac{1}{N(x_{i+1}(\delta)-x_i(\delta))}\int_{x_i(\delta)}^{x_{i+1}(\delta)}\varphi(x) dx\\
    & = \sum_{i\in \mathcal{S}^N}\frac{1}{N(x_{i+1}(\delta)-x_i(\delta))}\int_{x_i(\delta)}^{x_{i+1}(\delta)}\varphi(x) dx \\
    &\qquad + \sum_{i\in \{0,\ldots,N-1\}\setminus\mathcal{S}^N}\frac{1}{N(x_{i+1}(\delta)-x_i(\delta))}\int_{x_i(\delta)}^{x_{i+1}(\delta)}\varphi(x) dx,
\end{align*}
and, as $\delta\searrow0$, the right-hand side tends to
\begin{align*}
    & \frac{1}{N} \sum_{i\in \mathcal{S}^N}\varphi(\overline{x}_i) +  \sum_{i\in \{0,\ldots,N-1\}\setminus\mathcal{S}^N}\overline{R}_i\int_{\overline{x}_i}^{\overline{x}_{i+1}}\varphi(x) dx = \int_\R \varphi(x) d\overline{\rho}^N(x)\,.
\end{align*}
The convergence of $\rho^N(\cdot,\delta)$ to $\overline{\rho}^N$ in $2$-Wasserstein, as $\delta\searrow0$, immediately implies that the right-hand side of \eqref{eq:almost-weak-form2} converges as $\delta\searrow0$ to
\[
    -\int_\R \varphi(x,0)d \overline{\rho}^N(x) + \mathcal O(T N^{-1} \varepsilon^{-2}),
\]
due to the uniform continuity of $\varphi$. Then, we observe that 
\begin{align}
    \label{eq:delta_estimate}
    \int_0^\delta\int_\R \left|[\partial_t \varphi(x,t)  -   \Weps' \ast \rho^N(x, t) \partial_x \varphi(x,t) \right|\rho^N(x,t) dx dt \leq \delta \left(\norm{\partial_t \varphi}_{L^\infty}+\norm{\Weps'}_{L^\infty} \norm{\partial_x \varphi}_{L^\infty}\right).
\end{align}
Combining the above information we get, as $\delta\searrow0$ in \eqref{eq:almost-weak-form2}
\begin{align}
    \label{eq:almost-weak-form3}
    \int_0^T\int_\R \big[\partial_t \varphi(x,t)  -   \Weps' \ast \rho^N(x, t) \partial_x \varphi(x,t) \big] \rho^N(x,t)  dx dt &=  - \int_\R \varphi(x,0) d\overline{\rho}^N(x) dx\nonumber\\
    &\qquad + \mathcal O(T N^{-1} \varepsilon^{-2}).
\end{align}
We now want to pass to the limit $N\rightarrow+\infty$. Estimate \eqref{eq:discrete_smoothing} in Proposition \ref{prop:discrete_smoothing} implies that $\rho^N$ is uniformly bounded on every set of the form $\R\times [\sigma,T]$ for every fixed $\sigma>0$. Recalling Lemma \ref{lem:initial_convergence_measure} and using $\Weps' \in L^1(\R)$, arguing as in \eqref{eq:delta_estimate} we obtain the existence of a $\rho\in L^\infty_{\mathrm{loc}}((0,T];\,L^\infty(\R))$ such that for all $\sigma>0$
\begin{align*}
    & \int_0^T\int_\R \big[\partial_t \varphi(x,t)  -   \Weps' \ast \rho(x, t) \partial_x \varphi(x,t) \big] \rho(x,t)  dx dt =  - \int_\R \varphi(x,0) d\overline{\rho}(x) dx+ \mathcal{O}(\sigma).
\end{align*}
Since $\sigma>0$ was arbitrary, we may infer that $\rho$ is a weak solution to \eqref{eq:main_intro} in the sense of Definition \ref{def:weak_nonlocal}. Once again as in the last part of the proof of Theorem \ref{thm:convergence1}, $\rho$ is unique thanks to Theorem \ref{thm:existence_and_uniqueness_GF}, which proves the assertion.
\end{proof}

\section{The joint many-particle / small-interaction-range limit}\label{sec:joint_limit}

This section is devoted to proving the \textit{joint many particle / small interaction range limit} for initial datum in $L^1(\R)\cap L^\infty (\R)$.
The limit procedure requires specific constraints on the speed at which $\varepsilon\searrow 0$ and $N\rightarrow +\infty$. We provide them here in detail. Throughout the whole section we assume
\begin{itemize}
    \item [(J1)] $N\rightarrow+\infty$;
    \item [(J2)] $\varepsilon=\varepsilon_N$, where $\varepsilon_N>0$ is a sequence converging to zero;
    \item [(J3)] $\frac{1}{N\varepsilon_N^3}$ is uniformly bounded as $N\rightarrow+\infty$.
\end{itemize}

We provide here the precise statement of our result, which is the main one of this paper. 

\begin{thm}\label{thm:main}
   Let $\overline{\rho}\in \mathcal{P}_2 (\R)\cap L^\infty (\R)$. Assume that (J1)-(J2)-(J3) are satisfied. Then, the piecewise constant reconstruction of the density \eqref{eq:discrete_density} converges to the solution of the porous medium equation \eqref{eq:QPME}, in the sense of Definition \ref{def:weak_PME}, strongly in $L^p$, with $p\in [1,\infty]$.
\end{thm}
The dependence of various quantities on $\varepsilon$ is crucial in this context. The piecewise constant reconstruction of the density \eqref{eq:discrete_density} shall be denoted by $\rho^N_{\varepsilon_N}$ throughout. However, for simplicity in the notation, we will use the notation $\rhoen$ instead of $\rho^{N}_{\varepsilon_N}$. 
Moreover, the quantities $x_i$, $d_i$, and $R_i$ will still be written by omitting their dependence on $N$ and $\varepsilon_N$. We shall prove Theorem \ref{thm:main} following the strategy:
\begin{itemize}
    \item As a first step, we use Corollary \ref{cor:rho_Lp} to control $\rhoen$ uniformly in $L^2_{x,t}$. Up to a subsequence, we denote the weak limit by $\rho_0$.
    \item Using the elliptic law \eqref{eq:morse_elliptic}, we prove that $\Weps\ast\rho_0^N$ converges to $\rho_0$ weakly in $L^2 ([0,T]; L^2 (\R))$.
    \item We then show that $\Weps' \ast \rhoen$ is uniformly bounded in $L^2 (\R\times[0,T])$, by estimating the first moment and the functional $\int_\R \rhoen \log \rhoen\,dx$.
    \item We prove that $\Weps \ast \partial_t \rhoen$ is uniformly bounded in $L^2 ([0,T]; H^{-1}(\R))$. Therefore, by Aubin-Lions Lemma we have that $\Weps \ast \rhoen$ is strongly compact in $L^2 (\R \times (0,T))$, and we show that its strong limit coincides with $\rho_0$.
    \item A triangulation technique shows that the limit $\rho_0$ is attained \emph{strongly} in $L^2_{\mathrm{loc}}$, and this allows to pass to the limit and to obtain the concept of weak solution for \eqref{eq:PME}. 
\end{itemize}

In all the results below, the \emph{limit} is taken as $N\rightarrow+\infty$, which implies $\varepsilon\searrow 0$ due to (J2) above. Moreover, the expression \emph{uniformly bounded} refers to quantities that are uniformly bounded with respect to $N$ and $\varepsilon$.

\begin{lem}[Weak limit of $\rhoen$] \label{lem:weak_lim}
    Let $T\geq 0$. Assume $\bar{\rho}\in \mptr\cap L^\infty(\R)$. Then, the sequence $\rho_\varepsilon^N (x,t)$ admits, up to a non-relabelled subsequence, a weak limit $\rho_0$ in $L^2 (\R \times [0,T])$.
\end{lem}
\begin{proof}
    Lemma \ref{lem:initial_control} provides a uniform control of $\norm{\rhoen(\cdot,0)}_{L^p}$ for all $p\in [1,+\infty]$. Hence, Corollary \ref{cor:rho_Lp} shows that $\rhoen(\cdot,t)$ is uniformly bounded in all $L^p(\R)$ for all $p\in [1,+\infty]$ and for all $t\in [0,T]$, and consequently in $L^2(\R\times [0,T])$. The assertion follows from Kakutani's Theorem.
\end{proof}

\begin{lem}[Weak limit of $\Weps\ast \rhoen$] \label{lem:lim_conv} 
    Let $T\geq 0$. Assume $\bar{\rho}\in \mptr\cap L^\infty(\R)$. Then, 
    \[
    \Weps\ast \rhoen \rightharpoonup \rho_0,
    \]
    weakly in $L^2 ([0,T];L^2(\R))$.
\end{lem}
\begin{proof}
    Consider the elliptic law \eqref{eq:morse_elliptic} and let $\phi \in C_c^\infty (\R \times (0,T))$ be a test function. Then
    \[
    -\varepsilon^2 \int_0^T\int_\R \phi \Weps'' \ast \rhoen \, dx\,dt+ \int_0^T\int_\R \phi \Weps \ast \rhoen \,dx\,dt = \int_0^T\int_\R \phi \rhoen\,dx\,dt,
    \]
    that is, integrating by parts,
    \begin{equation}
    \label{eq:ell_conv}
    - \varepsilon^2\int_0^T \int_\R \phi_{xx} \Weps \ast \rhoen \, dx\,dt+ \int_0^T\int_\R \phi \Weps \ast \rhoen \,dx\,dt = \int_0^T \int_\R \phi \rhoen\,dx\,dt.
    \end{equation}
    The first term in \eqref{eq:ell_conv} vanishes as $\varepsilon$ vanishes due to Young's inequality. More precisely,
    \[\left|
    -\varepsilon^2 \int_0^T \int_\R \phi_{xx} \Weps \ast \rhoen \,dx\,dt \right|\leq  \varepsilon ^2 \norm{\phi _{xx}}_{L^\infty} \norm{\Weps}_{L^1} \norm{\rhoen}_{L^1} = \varepsilon ^2 \norm{\phi _{xx}}_{L^\infty}\rightarrow 0\,.
    \]
    Young's inequality also implies that the convolution $\Weps \ast \rhoen$ is uniformly bounded in $L^2 ([0,T]; L^2 (\R))$, therefore there exists $\chi \in L^2 ([0,T]; L^2 (\R))$ such that
    \[
    \int_0^T \int_\R \phi \Weps \ast \rhoen\,dx\,dt \to \int_0^T \int_\R \phi \chi \,dx\,dt \,.
    \]
    Concerning the last term in \eqref{eq:ell_conv}, recalling the definition of $\rho_0$, we get
    \[
    \int_0^T \int_\R \phi \rhoen \,dx\,dt \to \int_0^T \int_\R \phi \rho_0 \,dx\,dt\,.
    \]
    Hence, $\chi = \rho_0$, which shows that the assertion holds true.
\end{proof}

\begin{lem}\label{lem:energy}
    Let $\overline{\rho} \in L^\infty (\R)$. Assume (J1)-(J2)-(J3) are satisfied. Then, there exists a constant $C$ depending only on the $L^\infty$ norm of $\overline{\rho}$ such that
    \begin{align}\label{eq:energy}
    \begin{aligned}
         & \frac{1}{2}\int_\R \rhoen (x,t) \Weps \ast \rhoen (x,t)\,dx +\int_0^t\int_\R \rhoen (x,s) ( \Weps ' \ast \rhoen (x,s))^2 \,dx ds \\
         & \ \ \leq Ct +  \frac{1}{2}\int_\R \rhoen (x,0) \Weps \ast \rhoen (x,0)\,dx\,.
         \end{aligned}
    \end{align}
\end{lem}
\begin{proof}
We compute
    \begin{align*}
        & \frac{1}{2} \frac{d}{dt} \int_\R \rhoen (x,t) \Weps \ast \rhoen (x,t)\,dx \\
        & = \frac{1}{2}\frac{d}{dt} \sum_{i=0}^{N-1} \int_{x_i}^{x_{i+1}}  \Weps\ast \rhoen(x,t) \,dx R_i \\
        &= \frac{1}{2} \frac{d}{dt} \sum_{i=0}^{N-1} \sum_{k=0}^{N-1} \int_{x_i}^{x_{i+1}} \int_{x_k}^{x_{k+1}} \Weps (x-y)\,dy\,dx R_i R_k \\
        & = \frac{1}{2} \sum_{i=0}^{N-1} \sum_{k=0}^{N-1} \bigg[ \int_{x_i}^{x_{i+1}} \int_{x_k}^{x_{k+1}} \Weps (x-y) \,dy\,dx \,(\dot{R}_i R_k+ R_i \dot{R}_k ) \\
        & \qquad \qquad \qquad + R_i R_k \bigg( \dot{x}_{i+1} \int_{x_k}^{x_{k+1}} \Weps (x_{i+1} -y)\,dy - \dot{x}_i \int_{x_k}^{x_{k+1}} \Weps (x_i -y)\,dy \\
        & \qquad \qquad \qquad \qquad \qquad + \dot{x}_{k+1}  \int_{x_i}^{x_{i+1}}  \Weps (x-x_{k+1} ) \,dx - \dot{x}_k \int_{x_i}^{x_{i+1}}  \Weps (x-x_k) \,dx \bigg) \bigg] \\
        & = \sum_{i=0}^{N-1} \sum_{k=0}^{N-1} \bigg[ \int_{x_i} ^{x_{i+1}} \int_{x_k}^{x_{k+1}} \Weps (x-y) \,dx\,dy \dot{R}_i R_k \\
        & \qquad + R_iR_k \bigg( \dot{x}_{i+1} \int_{x_k}^{x_{k+1}} \Weps (x_{i+1} - y)\,dy - \dot{x}_i \int_{x_k}^{x_{k+1}} \Weps (x_i -y)\,dy \bigg) \bigg],
    \end{align*}
    where the last equality follows since the kernel $\Weps$ is even and we are summing $i$ and $k$ over the same set of indices. Thus, using \eqref{eq:scheme_rewritten} and \eqref{eq:R_i_dot} we get
    \begin{align*}
        & \frac{1}{2} \frac{d}{dt} \int_\R \rhoen (x,t) \Weps \ast \rhoen (x,t)\,dx \\
        & = \sum_{i=0}^{N-1} \sum_{k=0}^{N-1} \frac{R_i}{d_i} R_k \int_{x_i}^{x_{i+1}} \Weps'' \ast \rhoen(z,t) \,dz \int_{x_i}^{x_{i+1}} \int_{x_k}^{x_{k+1}} \Weps (x-y) \,dy\,dx \\
        & \qquad - R_i R_k \int_{x_i}^{x_{i+1}} \frac{\partial}{\partial z} \bigg( \Weps' \ast \rhoen (z,t) \int_{x_k}^{x_{k+1}} \Weps (z-y)\,dy \bigg) \,dz \\ 
        & = \sum_{i=0}^{N-1} \sum_{k=0}^{N-1} R_i R_k \int_{x_i}^{x_{i+1}} \Weps '' \ast \rhoen (z,t) \bigg[ \int_{x_i}^{x_{i+1}} \frac{1}{d_i} \int_{x_k}^{x_{k+1}} \Weps (x-y) \,dy\,dx \\
        & \hspace{9cm} - \int_{x_k}^{x_{k+1}} \Weps (z-y)\,dy \bigg] \,dz \\
        & \quad - \sum_{i=0}^{N-1} \sum_{k=0}^{N-1} R_i R_k \int_{x_i}^{x_{i+1}} \int_{x_k}^{x_{k+1}} \Weps' \ast \rhoen (z,t) \Weps ' (z-y) \,dy\,dz.
        \end{align*}
    
    We now estimate the first term on the right-hand side. We define $\Phi(x,t)$ such that
        \[
            \partial_x\Phi (x,t) = \sum_{k=0}^{N-1} R_k \int_{x_k}^{x_{k+1}} \Weps (x-y)\,dy = \Weps \ast \rhoen (x,t),
        \]
        and we obtain
        \begin{align*}
        & \sum_{i=0}^{N-1} \sum_{k=0}^{N-1} R_i R_k \int_{x_i}^{x_{i+1}} \Weps '' \ast \rhoen (z,t) \bigg[ \int_{x_i}^{x_{i+1}} \frac{1}{d_i} \int_{x_k}^{x_{k+1}} \Weps (x-y) \,dy\,dx - \int_{x_k}^{x_{k+1}} \Weps (z-y)\,dy \bigg] \,dz \\
        & = \sum_{i=0}^{N-1} R_i \int_{x_i}^{x_{i+1}} \Weps'' \ast \rhoen (z,t) \bigg[ \sum_{k=0}^{N-1} R_k \bigg( \int_{x_i} ^{x_{i+1}} \frac{1}{d_i} \int_{x_k}^{x_{k+1}} \Weps (z-y) \,dy\,dx \\
        & \hspace{9cm} - \int_{x_k}^{x_{k+1}} \Weps (z-y)\,dy \bigg) \bigg]\,dz \\
        & = \sum_{i=0}^{N-1} R_i \int_{x_i}^{x_{i+1}} \Weps'' \ast \rhoen (z,t) \bigg[ \frac{\Phi(x_{i+1},t) - \Phi (x_i,t)}{x_{i+1}-x_i} - \partial_z \Phi (z,t)  \bigg] \,dz.
    \end{align*}
   Combing the above computations we get
    \begin{align*}
        & \frac{1}{2} \frac{d}{dt} \int_\R \rhoen (x,t) \Weps \ast \rhoen (x,t)\,dx \\ 
        & = \sum_{i=0}^{N-1} R_i \int_{x_i}^{x_{i+1}} \Weps'' \ast \rhoen (z,t) \bigg[ \frac{\Phi(x_{i+1},t) - \Phi (x_i,t)}{x_{i+1}-x_i} - \partial_z \Phi (z,t)  \bigg] \,dz\\
        & \qquad - \sum_{i=0}^{N-1} \sum_{k=0}^{N-1} R_i R_k \int_{x_i}^{x_{i+1}} \int_{x_k}^{x_{k+1}} \Weps' \ast \rhoen (z,t) \Weps ' (z-y) \,dy\,dz \\
        & = \sum_{i=0}^{N-1} \int_{x_i}^{x_{i+1}} \rhoen(z,t) \Weps'' \ast \rhoen (z,t) \bigg[ \frac{\Phi(x_{i+1},t) - \Phi (x_i,t)}{x_{i+1}-x_i} - \partial_z \Phi (z,t)  \bigg] \,dz \\
        & \quad - \int_\R \rhoen (z,t) ( \Weps ' \ast \rhoen (z,t))^2 \,dz.
    \end{align*}
    We set
    \[
        T_\varepsilon \coloneqq \sum_{i=0}^{N-1} \int_{x_i}^{x_{i+1}} \rhoen(z,t) \Weps'' \ast \rhoen (z,t) \bigg[ \frac{\Phi(x_{i+1},t) - \Phi (x_i,t)}{x_{i+1}-x_i} - \partial_z \Phi (z,t)  \bigg] \,dz.
    \]
    and, proceeding as in the proof of Proposition \ref{prop:main_moment_estimate}, a Taylor expansion shows
    \begin{align*}
        & \bigg| \frac{\Phi (x_{i+1},t) - \Phi (x_i,t)}{x_{i+1} - x_i} - \partial_z \Phi (z,t) \bigg| \leq \frac{1}{2} \norm{\partial^2_{xx}\Phi }_\infty (x_{i+1}- x_i ) \\
        & \ = \frac{1}{2} \norm {\Weps'\ast\rhoen}_\infty (x_{i+1}- x_i )\leq  \frac{1}{2}\norm{\Weps'}_{L^1(\R)}\norm{\rhoen}_{L^\infty(\R)}d_i.
    \end{align*}
    Therefore, we control $T_\varepsilon$ by
    \begin{align*}
        \abs{T_\varepsilon} & \leq \frac{C}{2} \norm{\Weps'}_{L^1(\R)}\sum_{i=0}^{N-1} d_i \int_{x_i}^{x_{i+1}} \rhoen (x) \abs{\Weps'' \ast \rhoen (x)} \,dx \\
        & \leq \frac{C}{2N\varepsilon}\int_\R \abs{\Weps'' \ast \rhoen (x)} \,dx, 
    \end{align*}
    where $C$ depends on the $L^\infty$ norm or the initial condition $\overline{\rho}$. Hence, we estimate the last integral term as in the proof of Proposition \ref{prop:main_moment_estimate} and get
    \begin{align*}
        & \abs{T_\varepsilon}\leq \frac{C}{N\varepsilon^3}\,,
    \end{align*}
which proves the assertion due to assumption (J3).
\end{proof}

\begin{lem}[First moment estimate]\label{lem:first_moment}
  Let $\bar{\rho}\in \mptr\cap L^\infty(\R)$. Assume (J1)-(J2)-(J3) are satisfied. Then, there exists a constant $C$ such that, for all $t\geq 0$,
  \begin{equation}\label{eq:first_moment}
      \int_\R \abs{x}\rhoen(x,t) \,dx \leq \int_\R \abs{x} \rhoen(x,0) \, dx + \int_0^t \int_\R\rhoen(x,t) \abs{ \Weps'\ast\rhoen(x,s)}\,dx\, ds + Ct.
  \end{equation}
\end{lem}
\begin{proof}
    Taking $\phi(x)=\abs{x}$ in Proposition \ref{prop:main_moment_estimate}, we obtain
\begin{align*}
    & \frac{d}{dt}\int_\R \abs{x}\rhoen(x,t) dx = -\int_\R \rhoen(x,t) \sign(x)\Weps'\ast\rhoen(x,t)dx + S_\varepsilon,
\end{align*}
where
\[
\abs{S_\varepsilon}\leq \frac{[\phi]_{\mathrm{Lip}}}{N\varepsilon^2} = \frac{1}{N\varepsilon^2}\,.
\]
Hence, the statement easily follows from (J3) by taking the absolute value into the dissipation term. Note that, while $\phi$ does not belong to the class of test functions considered in Proposition \ref{prop:main_moment_estimate}, we can extend the class to Lipschitz continuous functions via a simple approximation procedure, and we omit the details.
\end{proof}

We now combine the estimates in Lemmas \ref{lem:energy}, \ref{lem:first_moment} and Corollary \ref{cor:rho_logrho} as follows.

\begin{lem}\label{lem:combined_log_energy_moment}
Let $\bar{\rho}\in \mptr\cap L^\infty(\R)$ and $T \geq 0$. Assume that $\int_\R \rhoen (x,0) \log \rhoen (x,0)\,dx < + \infty$. Then, there exists a constant $C>0$ only depending on the initial datum $\overline{\rho}$ and $T$ such that, for all $t \in [0,T]$,
    \begin{equation}\label{eq:log_dissipation}
         -\int_0^t\int_\R\rhoen(x,\tau)\Weps''\ast\rhoen(x,\tau)\,dx \, d\tau \leq C.
    \end{equation}
\end{lem}

\begin{proof}
From Lemma \ref{lem:estimates_1}, we know that $\rhoen(\cdot,t)$ is uniformly bounded in $L^2(\R)$ with respect to $N$ and $\varepsilon$ in terms of the $L^2$ norm of $\overline{\rho}$. By Young's inequality for convolutions we obtain that there exists a constant $C$ only depending on $\overline{\rho}$ such that
\[
\int_0^t\int_\R \rhoen(x,s) (\Weps'\ast\rhoen(x,s))^2\, dx\, ds \leq C(t+1)\,.
\]
By using Lemma \ref{lem:first_moment}, H\"older's inequality and the $L^\infty$ control on $\rhoen$ in Corollary \ref{cor:rho_Lp}, we obtain
\[
\int_\R \abs{x}\rhoen(x,t) \,dx  \leq C,
\]
for a suitable constant depending on $t$ and on $\overline{\rho}$.
Now, for an arbitrary constant $c>0$, we compute
    \begin{align*}
        & \int_\R \left(\rhoen(x,t) \log \rhoen(x,t) + \rhoen (x,t) c|x|\right) dx\\
        &\qquad = \int_\R \rhoen(x,t) (\log \rhoen(x,t) + c|x| )  \,dx \\
        &\qquad = \int_\R \rhoen(x,t) (\log \rhoen(x,t) - \log e^{-c |x|} ) \,dx\\
        &\qquad = \int_\R \frac{\rhoen(x,t)}{e^{-c|x|}} \log \frac{\rhoen(x,t)}{e^{-c|x| }} e^{-c|x|} \,dx.
    \end{align*}
We set $h(x,t)= \rhoen(x,t) / e^{-c|x|}$ and we choose $c$ such that 
\[
   \int_\R e^{- c|x|}\,dx = 1 = \int_\R \rhoen(x,t)\,dx,
\]
for all $t$. Therefore,
    \begin{align*}
        & \int_\R \left(\rhoen(x,t) \log \rhoen(x,t) + \rhoen (x,t) c|x|\right)\,dx
        \\ 
        &\qquad  = \int_\R h(x,t)\log h(x,t) e^{-c|x|}\,dx\\
        &\qquad = \int_\R e^{-c|x|} \psi(h(x,t))\,dx,
    \end{align*}
    with
    \[\psi(h)=h\log h - h +1\,.\]
    We observe that $\psi(h)\geq 0$ for all $h\geq 0$. The uniform estimate on the first moment and the result in Corollary \ref{cor:rho_logrho} imply
    \begin{align*}
        & \int_\R e^{-c|x|} \psi(\rhoen(x,t))\,dx -\int_0^t\int_\R\rhoen(x,\tau)\Weps''\ast\rhoen(x,\tau)dxd\tau \leq C,
    \end{align*}
for some $C$ depending only on $t$ and on the initial condition. The assertion follows from $\psi\geq 0$.
\end{proof}

We use the result in Lemma \ref{lem:combined_log_energy_moment} to control the $L^2$-norm of $\Weps'\ast\rhoen$.

\begin{prop} \label{prop:dissipation_estimate}
Assume $\bar{\rho}\in \mptr\cap L^\infty(\R)$. There exists a constant $C>0$ depending on $T$ and $\norm{\overline{\rho} }_{L^\infty(\R)}$ such that, for all $t\in [0,T]$,
\begin{equation}\label{eq:localised_2}
    \int_0^t \int_\R (\Weps'\ast\rhoen (x,s))^2 \,dx\,ds \leq  C\,.
\end{equation}
\end{prop}

\begin{proof}
Integration by parts yields
\begin{align}
    & \int_0^t \int_\R ( \Weps'\ast\rhoen(x,t))^2 dx ds  = - \int_0^t \int_\R \Weps \ast \rhoen (x,s) \Weps'' \ast \rhoen (x,s)\,dx\,ds ,\label{eq:localised_proof}
    \end{align}
where we have used that $\Weps\ast \rhoen$ decays to zero exponentially as $\abs{x}\rightarrow +\infty$, which justifies the integration by parts using that $\rhoen$ is compactly supported. 
Using the elliptic equation \eqref{eq:morse_elliptic}, we obtain
\begin{align*}
& \int_0^t \int_\R ( \Weps'\ast\rhoen(x,t))^2 \,dx\, ds \\
    & \ = - \int_0^t \int_\R \Weps \ast \rhoen (x,s) \Weps'' \ast \rhoen (x,s)\,dx\,ds  \\
    & \ = - \int_0^t \int_\R (\Weps \ast \rhoen (x,s) - \rhoen (x,s) ) \Weps'' \ast \rhoen (x,s) \,dx\,ds \\
    & \qquad - \int_0^t \int_\R \rhoen (x,s) \Weps'' \ast \rhoen (x,s)\,dx\,ds \\
    & \ = - \varepsilon^2 \int_0^t \int_\R ( \Weps'' \ast \rhoen (x,s) )^2\,dx\,ds - \int_0^t \int_\R \rhoen (x,s) \Weps'' \ast \rhoen (x,s)\,dx\,ds \\
    & \ \leq - \int_0^t \int_\R \rhoen (x,s) \Weps'' \ast \rhoen (x,s) \,dx\,ds,
\end{align*}
and the result follows from Lemma \ref{lem:combined_log_energy_moment}.
\end{proof}

\begin{lem}[$H^{-1}$ estimate of $\Weps\ast\partial_t \rho_\varepsilon^N$]\label{lem:negative_sobolev_1}
    Assume $\overline{\rho} \in \mathcal{P}_2 (\R) \cap L^\infty (\R)$. There exists a constant $C$ depending only on $\|\overline{\rho}\|_{L^\infty(\R)}$ such that
    \[
        \int_0^T \|\partial_t \Weps\ast\rhoen(\cdot,t)\|^2_{H^{-1}(\R)} dt \leq C,
    \]
    for all $N\in \mathbb{N}$ and for all $\varepsilon>0$ satisfying (J1)-(J2)-(J3).
\end{lem}

\begin{proof}
Let $T\geq 0$. Let $\phi\in H^1(\R)$ be a test function.
A computation similar to the one in the proof of Proposition \ref{prop:main_moment_estimate} shows
\begin{align}
    & \int_0^T \left|\langle \Weps \ast \partial_t \rhoen,\phi(\cdot)\rangle\right|^2 \, dt = \int_0^T \left| \int_\R \Weps\ast\partial_t\rhoen(x,t)\phi(x) \,dx\right|^2 \,dt  =\int_0^T \left|\int_\R \Weps\ast\phi(x)\partial_t\rhoen(x,t) \,dx\right|^2 \,dt \nonumber\\
    & \ \leq 2 \int_0^T\left|\int_{\R}\rhoen(x,t) \Weps'\ast\rhoen (x,t) \Weps\ast \phi'(x) \,dx \right|^2 \,dt + R(N,\varepsilon,T,\phi)\label{eq:negative_sobolev}
\end{align}
with
\[
R(N,\varepsilon,T,\phi)\leq 2 \frac{CT}{\varepsilon^4 N^2}\|\Weps\ast \phi'\|_{L^\infty(\R)}^2\leq \frac{CT}{\varepsilon^{5}N^2}\|\phi\|_{H^1(\R)}^2\,,
\]
where we have used $\|\Weps\|_{L^2(\R)}=\mathcal{O}\left(\varepsilon^{-1/2}\right)$. The first term on the right-hand side of \eqref{eq:negative_sobolev} is controlled, via Hölder's inequality and Young's inequality for convolutions, by
\begin{align*}
    & \sup_{t\geq 0}\norm{\rhoen(\cdot,t)}_{L^\infty(\R)}\int_0^T\left(\norm{\Weps\ast\phi'}_{L^2(\R)}^2 \left\|\Weps'\ast \rhoen(x,t)\right\|_{L^2(\R)}^2  \right) \,dt\\
    & \ \leq C \norm{\phi}_{H^1(\R)}^2 \int_0^T\int_\R\left(\Weps'\ast \rhoen(x,t) \right)^2 dx dt,
\end{align*}
and the result follows from Proposition \ref{prop:dissipation_estimate} and the $L^\infty$ estimate in Corollary \ref{cor:rho_Lp}.
\end{proof}

In view of the above results, the family $\Weps\ast \rhoen$ is uniformly bounded in
\[L^2([0,T];\, H^1(\R))\,,\]
whereas $\partial_t \Weps\ast \rhoen$ is uniformly bounded in 
\[L^2([0,T];\, H^{-1}(\R))\,.\]
Hence, Aubin-Lions Lemma implies the existence of a $\chi \in L^2 ((0,T); L^2 (\R))$ such that
\[
\Weps\ast \rhoen \to \chi
\]
strongly in $L^2_{\mathrm{loc}} (\R \times(0,T))$ (and, up to a subsequence, a.e.\ on $\R \times (0,T)$) in the joint $\varepsilon$-$N$ limit. In the next proposition, we will identify $\chi$ as the weak limit $\rho_0$ provided in Lemma \ref{lem:weak_lim}.

\begin{prop}\label{prop:strong_lim}
Assume $\overline{\rho}\in \mathcal{P}_2 (\R) \cap L^\infty (\R)$ and (J1)-(J2)-(J3) are satisfied. Let $\rho_0$ be the weak limit of $\rhoen$. Then, the strong limit of $\Weps\ast\rhoen$ in $L^2_{\mathrm{loc}} (\R \times (0,T))$ is $\rho_0$.
\end{prop}
\begin{proof}
    Let $\varphi\in C_c^\infty (\R \times (0,T))$. Then
    \[
    \int_0^T \int_\R \Weps\ast \rhoen \varphi\,dx\,dt \to \int_0^T \int_\R \chi \varphi\,dx\,dt.
    \]
    We also have that
    \[
    \int_0^T \int_\R \Weps\ast \rhoen \varphi\,dx\,dt = \int_0^T \int_\R \rhoen \Weps\ast \varphi \,dx\,dt,
    \]
    thus, we get
    \[
    \int_0^T \int_\R \rhoen \Weps\ast \varphi \,dx\,dt \to \int_0^T \int_\R \chi \varphi\,dx\,dt.
    \]
    We recall that $\rhoen \rightharpoonup \rho_0$ in $L^2((0,T);L^2(\R))$. Moreover, $\Weps\ast \varphi \to\varphi$ strongly. To see this, we observe 
    \[\|\Weps\ast \varphi-\varphi\|_{L^2} = \varepsilon^2\|\Weps''\ast\varphi\|_{L^2} = \varepsilon^2 \|\Weps\ast \varphi''\|_{L^2}\leq \varepsilon^2 \|\Weps\|_{L^1}\|\varphi''\|_{L^2}\rightarrow 0\,.\]
Hence,
    \[
    \int_0^T \int_\R \rhoen \Weps\ast \varphi\,dx\,dt \to \int_0^T \int_\R \rho_0 \varphi\,dx\,dt
    \]
    and $\chi=\rho_0$.
\end{proof}

We are now in the position to prove the strong $L^2$ compactness of $\rhoen$.

\begin{prop}
    Assume $\overline{\rho} \in \mathcal{P}_2 (\R) \cap L^\infty (\R)$ and (J1)-(J2)-(J3) are satisfied.
    It holds that $\rhoen \to \rho_0$ in $L^2_{\mathrm{loc}} (\R\times(0,T))$ in the joint limit $\varepsilon \to 0$ and $N \to +\infty$.
\end{prop}
\begin{proof}
    Using the elliptic law \eqref{eq:morse_elliptic} and integration by parts, we get
    \begin{align*}
        \int_0^T \int_\R \abs{\rhoen - \Weps \ast \rhoen}^2\,dx\,dt & = \int_0^T \int_\R  (\rhoen - \Weps\ast\rhoen) (\rhoen - \Weps \ast \rhoen)\,dx\,dt \\
        & = - \varepsilon^2 \int_0^T \int_\R  \Weps'' \ast \rhoen (\rhoen - \Weps \ast \rhoen)\,dx\,dt \\
        & = - \varepsilon^2 \int_0^T \int_\R  \rhoen \Weps '' \ast \rhoen\,dx\,dt + \varepsilon^2 \int_0^T \int_\R \Weps'' \ast\rhoen \Weps \ast \rhoen \,dx\,dt \\
        & = - \varepsilon^2 \int_0^T \int_\R \rhoen \Weps '' \ast \rhoen\,dx\,dt - \varepsilon^2\int_0^T \int_\R (\Weps' \ast \rhoen )^2\,dx\,dt\\
        & \leq - \varepsilon^2 \int_0^T \int_\R \rhoen \Weps '' \ast \rhoen\,dx\,dt\,.
    \end{align*}
    From Lemma \ref{lem:combined_log_energy_moment} we have that the first term on the r.h.s.\ vanishes, thus
    \[ \int_0^T \int_\R \abs{\rhoen - \Weps \ast \rhoen }^2\,dx\,dt \leq \mathcal{O}(\varepsilon^2) \,.\]
    Therefore, on a compact set $K\subset \R$, we have
    \begin{align*}
        \norm{\rhoen - \rho_0}_{L^2((0,T); L^2(K))} & \leq \norm{\rhoen - \Weps\ast \rhoen}_{L^2((0,T); L^2(K))} + \norm{\Weps \ast \rhoen - \rho_0}_{L^2((0,T); L^2(K))} \\
        & \leq \mathcal{O}(\varepsilon) + \norm{\Weps \ast \rhoen - \rho_0}_{L^2((0,T); L^2(K))}.
    \end{align*}
    By using Proposition \ref{prop:strong_lim}, we have the statement.
\end{proof}

We now establish that the limit $\rho_0$ is a weak solution to the porous medium equation, thus completing the proof of Theorem \ref{thm:main}.

\begin{proof}[Proof of Theorem \ref{thm:main}]
    Let $\varphi\in C_c^\infty (\R \times (0,T))$ be a test function. The sequence $\rhoen$ fulfils, see Definition \ref{def:weak_nonlocal},
    \begin{equation}
    \label{eq:proof_weak_nl}
    \int_0^T \int_\R \partial_t \varphi \rhoen \,dx\,dt + \int_0^T \int_\R \partial_x \varphi \rhoen \Weps' \ast \rhoen \,dx\,dt=\mathcal{O}\left(\frac{1}{\varepsilon N^2}\right).
    \end{equation}
    The first term in \eqref{eq:proof_weak_nl} converges to $\int_0^T \int_\R \partial_t \varphi\rho_0\,dx\,dt$, due to the definition of $\rho_0$.
    Concerning the second term, we have that $\Weps'\ast \rhoen \rightharpoonup \partial_x \rho_0$ weakly in $L^2$. Indeed
    \[
    \int_0^T \int_\R \varphi \Weps'\ast \rhoen \,dx\,dt = - \int_0^T \int_\R \partial_x \varphi \Weps \ast \rhoen\,dx\,dt,
    \]
    and then we use Proposition \ref{prop:strong_lim}. Therefore, since $\rhoen \to \rho_0$ strongly in $L^2$, we obtain that the second term in \eqref{eq:proof_weak_nl} converges to 
    \[
    - \int_0^T \int_\R \partial_x \varphi \rho_0 \partial_x \rho_0\,dx\,dt \,.
    \]
 By Definition \ref{def:weak_PME}, we have the statement up to a subsequence. Classical results on the uniqueness of weak solutions to the porous medium equations, see \cite{vazquez_book_1}, implies the whole sequence converges strongly to $\rho_0$.
\end{proof}

\section*{Acknowledgments}
The research of MDF is supported by the Ministry of University and Research (MIUR), Italy under the grant PRIN 2020- Project N. 20204NT8W4, Nonlinear Evolutions PDEs, fluid
dynamics and transport equations: theoretical foundations and applications. The research of VI is supported by the Italian INdAM project N. E53C22001930001 “MMEAN-FIELDS”.
MDF and VI acknowledge support from the InterMaths Network, \url{www.intermaths.eu}. Part of this work was carried out during a visit of MDF and VI at the Institute of Scientific Computing of TU Dresden, which provided kind hospitality and support to this research. MS would like to express his gratitude for the hospitality of the University of L'Aquila where part of the discussions took place.

\appendix

\section{Some technical results}\label{appendix1}

In this section we collect some technical results which are used throughout the paper. We start by proving Theorem \ref{thm:continuity}.

\begin{proof}[Proof of Theorem \ref{thm:continuity}]
    Let $\overline{\rho}\in \mptr$. For $t\in [0,T]$ we define the \emph{flow map} $Y(\cdot,t):\R\rightarrow\R$ as the solution to the Cauchy problem
\begin{equation}\label{eq:flow_map}
     \begin{dcases}
         \partial_t Y(x_0,t) \!\!\!\! &= v(Y(x_0,t),t), \\
         \phantom{..} Y(x_0,0) &=x_0,
     \end{dcases}
\end{equation}
for all $x_0\in \R$. The classical ODE theory implies existence and uniqueness of $Y(x_0, \cdot)$ for all $x_0\in \R$ due to the Lipschitz assumption on $v$. We set

\begin{equation}\label{eq:CE1}
\rho(\cdot,t)=Y(\cdot,t)_{\#}\overline{\rho}
\end{equation}
  for $t\in [0,T]$. Classical results (see for simplicity the recent \cite[Theorem 4.4]{santambrogio_book}) imply that $\rho$ defined in \eqref{eq:CE1} is the unique solution to the continuity equation \eqref{eq:continuity} on $[0,T]$ with initial condition $\rho(\cdot,0)=\overline{\rho}$.
  Now, let us set
  \[\overline{X}=\mathcal{T}(\overline{\rho})\]
   and consider the unique solution $X(z,t)$ to the Cauchy problem \eqref{eq:continuity_quantile} with $X(\cdot,0)=\overline{X}$. By integrating \eqref{eq:continuity_quantile} on $[0,t]$ we get, for all $t\in [0,T]$,
\begin{equation}\label{eq:mild_quantile}
    X(z,t)=\overline{X}(z)+\int_{0}^t v(X(z,\tau),\tau)\,d\tau\,.   
   \end{equation}
On the other hand, by integrating the ODE in \eqref{eq:flow_map} on $[0,t]$ with $t\in[0,T]$ 
 we get
   \[Y(x_0,t)=x_0+ \int_{0}^t v(Y(x_0,\tau),\tau)\,d\tau\,.\]
By setting $x_0=\overline{X}(z)$, we get that $X(z,t)$ and $Y(\overline{X}(z),t)$ 
are mild solutions (with same initial condition) to an ODE that admits a unique solution, and are therefore equal on $(z,t)\in [0,1]\times [0,T]$. To prove the first implication of the theorem we only need to prove that $X(z,t)=\mathcal{T}(\rho(\cdot,t))(z)$. Now, \eqref{eq:CE1} implies for all $\phi:\R\rightarrow\R$ measurable
\[\int_\R \phi(x)\rho(x,t) \,dx = \int_\R \phi(Y(x_0,t))\overline{\rho}(x_0)\, dx_0 = \int_0^1 \phi(Y(\overline{X}(z),t))\,dz\]
which is equivalent to 
\[
\rho(\cdot,t)=Y(\overline{X}(\cdot),t)_\#\mathcal{L}^1_{[0,1]},
\]
which implies
\[X(\cdot,t)=Y(\overline{X}(\cdot),t) = \mathcal{T}(\rho(\cdot,t))\,.\]

Vice versa, let $X\in L^\infty([0,T]; L^2([0,1]))$ solve \eqref{eq:continuity_quantile} a.e.\ on $[0,1]\times [0,T]$ with initial condition $\overline{X}\in L^2([0,1])$. Hence, $X(\cdot,\cdot)$ solves \eqref{eq:mild_quantile} on $(z,t)\in [0,1]\times [0,T]$. Let $\overline{\rho}\in \mptr$ given by $\overline{\rho}=\mathcal{T}^{-1}(\overline{X})$, and $\overline{F}=F_{\overline{X}}$ as in \eqref{eq_F_X}, which implies $\overline{\rho}=\partial_x \overline{F}$ in the sense of distributions. We write \eqref{eq:mild_quantile} with $z=\overline{F}(x_0)$ for $x_0\in \R$:
\begin{equation}\label{eq:flow_map_modified}
    X(\overline{F}(x_0),t)=\overline{X}(\overline{F}(x_0)) + \int_0^t v(X(\overline{F}(x_0),\tau),\tau) \,d\tau\,.
\end{equation}
Now, in general the identity $\overline{X}(\overline{F}(x_0))=x_0$ is not true, but one can easily verify that it is true for all $x_0\in \mathrm{supp}(\overline{\rho})$. Hence, for all $x_0\in \mathrm{supp}(\overline{\rho})$ we have
\[X(\overline{F}(x_0),t)=x_0+ \int_0^t v(X(\overline{F}(x_0),\tau),\tau) \,d\tau\,,\]
which means that $Y(x_0,t) \coloneqq X(\overline{F}(x_0),t)$ solves \eqref{eq:flow_map} for all $x_0\in \mathrm{supp}(\overline{\rho})$. This is enough to conclude, once again due to \cite[Theorem 4.4]{santambrogio_book}, that
\[\rho(\cdot,t)=Y(\cdot,t)_{\#}\overline{\rho} \]
is a weak solution to \eqref{eq:continuity} on $[0,T]$ with initial condition $\overline{\rho}$. Moreover, for all $\phi:\R\rightarrow\R$ measurable,
\begin{align*}
    & \int_\R \rho(x,t) \phi(x)\, dx = \int_\R\phi(Y(x_0,t))\overline{\rho}(x_0) \,dx_0 = \int_\R \phi(X(\overline{F}(x_0),t))\overline{\rho}(x_0)\,dx_0\\
    & \ = \int_0^1 \phi(X(z,t)) \,dz
\end{align*}
where we have used once again that $Y(x_0,t)$ and $X(\overline{F}(x_0),t)$ coincide on $\mathrm{supp}(\overline{\rho})$, the identity $\overline{\rho}=\overline{X}_{\#}\mathcal{L}^1_{[0,1]}$, and the identity $\overline{F}(\overline{X}(z))=z$, which is justified by the fact that $\overline{\rho}\in L^\infty(\R)$. Indeed, $\overline{\rho}\in L^\infty(\R)$ implies $\overline{X}=\mathcal{T}(\overline{\rho})$ is strictly increasing on $[0,1]$ and therefore $\overline{F}\circ\overline{X}(z)=z$ for all $z\in [0,1]$. The above identity shows that $\rho(\cdot,t)=\mathcal{T}^{-1}(X(\cdot,t))$ which proves the desired assertion.
\end{proof}

We next provide the proof of Lemma \ref{lem:N_and_R_coincide}, which was already contained in \cite{bonaschi_thesis} and \cite{BCDP}.

\begin{proof}[Proof of Lemma \ref{lem:N_and_R_coincide}]
Let $X\in \mathcal{K}$. Since $X$ is non-decreasing, recalling the definition of $\mathcal{N}_\varepsilon$, we compute
    \begin{align*}
        \mathcal{N}_\varepsilon [X] & =\frac{1}{2}\int_0^1\int_0^1 \Neps(X(z)-X(\zeta))\,d\zeta\, dz = \frac{1}{4\varepsilon}\int_0^1\int_0^1 \left(1-\frac{\abs{X(z)-X(\zeta)}}{\varepsilon}\right)\,d\zeta\, dz \\
        & = \frac{1}{4\varepsilon}\left(1-\frac{1}{\varepsilon}\int_0^1\int_0^1 \abs{X(z)-X(\zeta)}\,d\zeta \,dz \right)\\
        & = \frac{1}{4\varepsilon}\left[1-\frac{1}{\varepsilon}\left(\iint_{z>\zeta}(X(z)-X(\zeta)) \,d\zeta \,dz +\iint_{z<\zeta}(X(\zeta)-X(z)) \,d\zeta \,dz\right)\right]\\
        & = \frac{1}{4\varepsilon}\left[1-\frac{2}{\varepsilon}\iint_{z>\zeta}(X(z)-X(\zeta)) \,d\zeta \,dz\right]\\
        & =  \frac{1}{4\varepsilon}\left[1-\frac{2}{\varepsilon}\left(\int_0^1 \int_0^z X(z)\,d\zeta \,dz-\int_0^1\int_{\zeta}^1 X(\zeta)\,dz \,d\zeta\right)\right]\\
        & = \frac{1}{4\varepsilon}\left[1-\frac{2}{\varepsilon}\left(\int_0^1 z X(z) \,dz-\int_0^1(1-\zeta)X(\zeta)\, d\zeta\right)\right]\\
        & = \frac{1}{4\varepsilon}\left[1-\frac{2}{\varepsilon}\int_0^1(2z-1)X(z) \,dz\right]=\mathcal{R}_\varepsilon[X]\, ,
    \end{align*}
    thus the statement is proven.
\end{proof}


\end{document}